\pdfoutput=1
\documentclass[a4paper, DIV12, pointlessnumbers]{scrartcl}

\usepackage[utf8x]{inputenc}
\usepackage[T1]{fontenc}
\usepackage[english]{babel}
\usepackage{lmodern}

\usepackage{amsmath}
\usepackage{amssymb}
\usepackage{amstext}
\usepackage{amsfonts}
\usepackage{array}
\usepackage{theorem}

\usepackage[longnamesfirst]{natbib}

\numberwithin{equation}{section}

\theoremheaderfont{\sf\bfseries}
\newtheorem{definition}{Definition}[section]
\newtheorem{proposition}[definition]{Proposition}
\newtheorem{lemma}[definition]{Lemma}
\newtheorem{theorem}[definition]{Theorem}
\newtheorem{corollary}[definition]{Corollary}
\newtheorem{assumption}[definition]{Assumption}

\theorembodyfont{\rmfamily}
\newtheorem{remark}[definition]{Remark}
\newtheorem{exmpl}[definition]{Example}

\renewcommand{\star}{\ast}

\renewcommand{\epsilon}{\varepsilon}

\newcommand{\NN}{\mathbb{N}}

\newcommand{\RR}{\mathbb{R}}

\newcommand{\cF}{\mathcal{F}}

\newcommand{\N}{\mathcal{N}}
\newcommand{\cN}{\mathcal{N}}

\newcommand{\R}{\mathcal{R}}

\newcommand{\Z}{\mathcal{Z}}

\newcommand{\hf}{\widehat{f}}

\newcommand{\hk}{\widehat{k}}

\newcommand{\hw}{\widehat{w}}

\newcommand{\td}{\widetilde{d}}

\newcommand{\tr}{\widetilde{r}}

\newcommand{\tw}{\widetilde{w}}

\newcommand{\E}{\mathbf{E}}
\newcommand{\Ex}{\mathbf{E}}
\renewcommand{\P}{\mathbf{P}}

\newcommand{\var}{\operatorname{Var}}
\renewcommand{\var}{\operatorname{\mathbf{V}\mathrm{ar}}}

\newcommand{\argmin}{\operatorname*{argmin}}
\newcommand{\argmax}{\operatorname*{argmax}}

\renewcommand{\epsilon}{\varepsilon}

\DeclareMathOperator{\skalar}{\langle\cdot,\cdot\rangle}   
\newcommand{\skalarV}[1]{\langle #1\rangle}   
\DeclareMathOperator{\norm}{\lVert\cdot\rVert}   
\newcommand{\normV}[1]{\lVert#1\rVert}   

\newcommand{\I}[1]{\mathbf{1}_{#1}}
\newcommand{\1}{{\mathbf 1}}

\newcommand{\proof}{\noindent\textit{Proof.} }
\newcommand{\proofof}[1]{\noindent\textit{Proof of #1.}}

\newcommand{\qed}{\hfill$\Box$\par\vspace{1em}}
\newcommand{\eoe}{\hfill $\diamondsuit$}

\newcounter{nc}
\newcommand{\diag}{\operatorname{diag}}

\newcommand{\wrt}{with respect to\ }

\RequirePackage{color}
\definecolor{dgreen}{rgb}{0,0.5,0}
\definecolor{dblue}{rgb}{0,0,0.9}
\definecolor{dred}{rgb}{0.6,0.0,0.1}
\definecolor{dgold}{rgb}{0.9,0.4,0.0}
\definecolor{dvio}{rgb}{0.6,0.3,0.5}
\definecolor{gray}{rgb}{0.5,0.5,0.5}

\renewcommand{\leq}{\leqslant}
\renewcommand{\geq}{\geqslant}

\renewcommand{\phi}{\varphi}
\renewcommand{\epsilon}{\varepsilon}

\renewcommand{\P}{\mathbf{P}}

\renewcommand{\proof}{\noindent\textit{Proof.} }
\renewcommand{\qed}{\hfill $\Box$\par\vspace{1em}}

 \renewcommand{\N}{{\mathbb N}}

\newcommand{\cB}{\mathcal{B}}
\newcommand{\cS}{\mathcal{S}}

\newcommand{\cE}{{\cal E}}

\newcommand{\supt}[1]{\sup_{t\in\cB_{#1}}}

\newcommand{\htau}{\widehat{\tau}}

\newcommand{\pen}{\operatorname{pen}}
\newcommand{\hpen}{\widehat{\pen}}

\newcommand{\hPhi}{\widehat\Phi}
\newcommand{\tPhi}{\widetilde\Phi}
\newcommand{\ct}{\Upsilon}

\newcommand{\whk}{\widehat k}
\newcommand{\kstar}{k_n^*}

\newcommand{\Cy}{N}
\newcommand{\Ce}{M}
\newcommand{\hCy}{\widehat N}

\newcommand{\xdfw}{\gamma}
\newcommand{\xdfr}{r}
\newcommand{\edfq}{q}
\newcommand{\edfw}{\lambda}
\newcommand{\Edfw}{\Lambda}
\newcommand{\edfr}{d}

\newcommand{\xdf}{f}
\newcommand{\ydf}{g}
\newcommand{\edf}{\varphi}

\newcommand{\hydf}{\widehat\ydf}

\newcommand{\txdf}{\widetilde\xdf}

\newcommand{\qtII}{\skalarV{t,\hPhi_g - \tPhi_g}_\hw}

\newcommand{\bw}{\xdfw}
\newcommand{\br}{\xdfr}
\renewcommand{\td}{d}
\renewcommand{\tw}{\lambda}
\renewcommand{\hw}{\omega}

\newcommand{\hdelta}{\widehat{\delta}}
\newcommand{\hDelta}{\widehat{\Delta}}

\newcommand{\ii}{
\renewcommand{\theenumi}{(\roman{enumi})}
\renewcommand{\labelenumi}{\theenumi}
}

\newcommand{\arab}{
\renewcommand{\theenumi}{\arabic{enumi}.}
\renewcommand{\labelenumi}{\theenumi}
}

\renewcommand{\hk}{\whk}

\newcommand{\Fgr}{\cF_\xdfw^\xdfr}

\newcommand{\fou}[1]{[#1]}
\newcommand{\hfou}[1]{\widehat{[#1]}}

\newcommand{\tschebby}{Chebychev}

\newcommand{\fedf}{\fou{\edf}}

\newcommand{\hfedf}{\hfou{\edf}}

\begin{document}

\makeatletter%
\def\@fnsymbol#1{\ensuremath{\ifcase#1\or * \or \star \or 1 \or 2\or 3\or  , \or
g\or h\or i\else\@ctrerr\fi}}%
\makeatother%

\author{{\sc Jan Johannes}
\thanks{Institut de
    statistique, biostatistique et sciences actuarielles (ISBA), Voie du Roman Pays 20, 1348~Louvain-la-Neuve,
    Belgium, e-mail: {\{jan.johannes|maik.schwarz\}@uclouvain.be}}
  \and {\sc Maik
    Schwarz$\;^*$
}}

\title{Partially adaptive nonparametric\\ instrumental regression}

\date{Université catholique de Louvain}
\dedication{10 April 2012}
\maketitle

\begin{abstract} We consider the problem of estimating the structural
  function in nonparametric instrumental regression, where in the
  presence of an instrument $W$ a response~$Y$ is modeled in
  dependence of an endogenous explanatory variable $Z$.

 The proposed estimator is based on dimension reduction and additional
 thresholding. The minimax optimal rate of convergence of the
 estimator is derived assuming that the structural function belongs to
 some ellipsoids which are in a certain sense linked to the
 conditional expectation  of $Z$ given $W$.~We illustrate
 these results by considering classical smoothness assumptions.
 However, the proposed estimator requires an optimal choice of a
 dimension parameter depending on certain characteristics of the
 unknown structural function and the conditional expectation 
 of $Z$ given $W$, which are not known in practice.  The main issue
 addressed in our work is an adaptive choice of this dimension
 parameter using a model selection approach under the restriction that
 the conditional expectation of $Z$ given $W$ is smoothing in
 a certain sense.  In this situation we develop a penalized minimum
 contrast estimator with randomized penalty and collection of models.
 We show that this data-driven estimator can attain the lower risk
 bound up to a constant over a wide range of smoothness classes for
 the structural function.
\end{abstract}
{\footnotesize
\begin{tabbing}
\> Minimax theory,  Orthogonal series estimation, Model selection, Adaptive estimation.\\
\end{tabbing}
\textbf{Acknowledgment.} This work was supported by the IAP research network no.\ P6/03 of the
Belgian Government (Belgian Science Policy) and by the ``Fonds
Sp\'eciaux de Recherche'' from the Universit\'e catholique de Louvain.
}

\renewcommand{\asymp}{\sim}

 \newcommand{\cT}{{\cal T}}
 \newcommand{\cU}{{\cal U}}

 \newcommand{\cW}{{\cal W}}

 \newcommand{\uk}{\underline{k}}
 \newcommand{\hmu}{\widehat{\mu}}

\renewcommand{\xdfw}{\gamma}
\renewcommand{\xdfr}{r}
\renewcommand{\edfq}{q}
\renewcommand{\edfw}{\lambda}
\renewcommand{\Edfw}{\Lambda}
\renewcommand{\edfr}{d}

\renewcommand{\xdf}{f}
\renewcommand{\ydf}{g}
\renewcommand{\edf}{\varphi}

\newcommand{\ukstar}{\underline{\kstar}}
\newcommand{\dstar}{{R_n^*}}

 \newcommand{\wtk}{\widetilde k}

\newcommand{\rep}{h}
\newcommand{\ev}{\lambda}


\newcommand{\sol}{\varphi}
\newcommand{\Op}{T}
\newcommand{\gf}{g}
\newcommand{\basZ}{e}
\newcommand{\basW}{f}
\newcommand{\basTr}{\psi}

\newcommand{\solw}{\gamma}
\newcommand{\solr}{\rho}
\newcommand{\Opw}{\lambda}
\newcommand{\Opd}{d}
\newcommand{\OpD}{D}

\newcommand{\hsol}{\widehat\sol}
\newcommand{\tsol}{\widetilde\sol}
\newcommand{\hOp}{\widehat\Op}
\newcommand{\tOp}{\widetilde\Op}
\newcommand{\hgf}{\widehat\gf}
\newcommand{\tgf}{\widetilde\gf}

\newcommand{\cTdDw}{\cT_{\Opd,\OpD}^{\Opw}}
\newcommand{\cTdw}{\cT_{\Opd}^{\Opw}}

\renewcommand{\td}{d}
\renewcommand{\tw}{\lambda}

\newcommand{\an}{n^{1/3}}
\newcommand{\hnu}{\widehat\nu}

\newcommand{\hwtwD}{\Psi}


\newcommand{\vuz}{V_{U|Z}}
\newcommand{\Evuz}{\E[\vuz]}


\newpage
\section{Nonparametric instrumental regression}
\label{D:cha:NIV}

Nonparametric instrumental regression models have attracted
increasing attention in the statistic and econometric literature
\citep[e.g.][]{F03eswc,DFR02,NP03econometrica,HallHorowitz2007,BCK07econometrica}. In
many applications, the dependence of a response
$Y$ on the variation of an \textit{endogenous} vector $Z$  of explanatory
variables is characterized by
\begin{subequations}
  \begin{equation}
    \label{D:model:NP}
    Y = \varphi(Z) + U
  \end{equation}
  for some error term $U$. Endogenous means that $Z$ and $U$ are not
  stochastically mean-independent (i.e., $\E[U|Z]\neq0$). The
  nonparametric relationship is hence modeled by the regression
  function $\sol$, which is also called structural function. However, classical nonparametric regression methods cannot be
  applied directly in this context. The approach of instrumental
  regression copes with the mean dependence by introducing an additional
  vector of exogenous \textit{instruments}~$W$ such that
  \begin{equation}
    \label{D:model:NP2}
    \Ex [U | W] = 0.
  \end{equation}
\end{subequations}
 Typical examples of such settings are error-in-variable
models, simultaneous equations or treatment models with endogenous
selection. It is worth noting that in the presence of instrumental
variables, the model equations (\ref{D:model:NP}--\ref{D:model:NP2})
are the natural generalization of a standard parametric model
\citep[e.g.][]{A74je} to the nonparametric situation. This extension
has first been introduced by \cite{F03eswc} and
\cite{NP03econometrica}, while its identification has been studied
e.g.\ in \cite{CFR:07}, \cite{DFR02} and \cite{FJVB07}. Recent
applications and extensions of this approach include nonparametric
tests of exogeneity \citep{BH07res}, quantile regression models
\citep{HL07econometrica}, or semi-parametric modeling \citep{FJVB05},
for example.

 There is a vast literature on the nonparametric estimation
of the structural function~$\sol$ based on a sample of $(Y,Z,W)$. For
example, \cite{AC03econometrica}, \cite{BCK07econometrica} or
\cite{NP03econometrica} consider sieve minimum distance estimators,
while \cite{DFR02}, \cite{GS06} or \cite{FJVB07} consider penalized
least squares estimators. The optimal estimation in a minimax sense
has been worked on by \cite{HallHorowitz2005} and \cite{CR:10}. The
authors prove a lower bound for the mean integrated squared error
(MISE) and propose an estimator which can attain optimal rates. In the
present article, we extend this result by considering not only the
MISE of the estimation of $\sol$ but, more generally, a risk defined
\wrt a weighted norm. This allows us for example to
consider the estimation of the derivatives of $\sol$, too.

It is well known that all the resulting estimation procedures can
attain optimal rates only if certain smoothing parameters are chosen
in an appropriate way. In general, this choice requires knowledge of
characteristics of the structural function, such as the number of its
derivatives, which are not known in practice. Thus, an essential and
still open problem in this theoretical framework is the data driven
choice of smoothing parameters. In this paper, an adaptive method is
proposed which indeed does not depend on any properties of $\sol$,  though
yielding optimal rates. However, it still necessitates that
some characteristics of the underlying conditional expectation be known.

One objective in this article is the minimax optimal nonparametric
estimation of the structural function $\sol$ based on an iid.\ sample
of $(Y,Z,W)$ satisfying the model equations
(\ref{D:model:NP}--\ref{D:model:NP2}). 
Let us briefly sketch our estimation approach here.  For the moment
being, suppose that the structural function can be represented as
$\sol=\sum_{j=1}^k[\sol]_j \basZ_j$ using only $k$ pre-specified basis
functions $\basZ_1,\dotsc,\basZ_k$, and that only the coefficients
$[\sol]_j$ \wrt these functions are unknown.  In this situation,
taking the conditional expectation \wrt the instrument $W$ on both sides of~\eqref{D:model:NP}
yields a
multivariate linear conditional moment equation, that is, $\Ex[Y | W]
= \sum_{j=1}^k [\sol]_j\Ex[\basZ_j(Z)| W]$.  Solving this equation is
a classical textbook problem in econometrics \citep[cf.][]{PU99}. A
popular approach consists in replacing the conditional moment equation
by an unconditional one: given $k$ functions $\basW_1,\dotsc,\basW_k$,
one can consider~$k$ unconditional moment equations instead of the
multivariate conditional moment equation, that is, $\Ex[Y \basW_l(W)]
= \sum_{j=1}^k [\sol]_j\Ex[\basZ_j(Z)\basW_l(W)]$,
$l=1,\dotsc,k$. Notice that once the functions $\{\basW_l\}_{l=1}^k$
are chosen, all the unknown quantities in the unconditional moment
equations can be estimated by simply substituting empirical versions
for the theoretical expectation. Moreover, a least squares solution of
the estimated equation leads to a consistent and asymptotically
normally distributed estimator of the parameter vector
$([\sol]_j)_{j=1}^k$ under mild assumptions.  The choice of the
functions $\{\basW_l\}_{l=1}^k$ directly influences the asymptotic
variance of the estimator and thus the question of optimal instruments
minimizing the asymptotic variance arises
\citep[cf.][]{N90econometrica}. One advantage of this approach is that
the estimator is easily computable. However, in many situations an
infinite number of functions $\{\basZ_j\}_{j\geq1}$ and associated
coefficients $([\sol]_j)_{j\geq1}$ is needed to represent the
structural function $\sol$.  The choice of the basis functions
$\{\basZ_j\}_{j\geq1}$ reflects a priori information about the
structural function~$\sol$, such as smoothness. Considering an
infinite number of functions $\{\basW_l\}_{l\geq1}$, we could still
consider the finite dimensional least squares estimator described
above for each $k\geq1$.

Notice that the dimension $k$ plays the role of a smoothing parameter
and one might expect that the estimator of the structural function
$\sol$ is consistent as $k$ tends to infinity at a suitable
rate. Unfortunately, this is not true in general. Let
$\sol_k:=\sum_{j=1}^k [\sol_k]_j \basZ_j$ denote a least squares
solution of the reduced unconditional moment equations. This means
that the vector of coefficients $([\sol_k]_j)_{j=1}^k$ minimizes the
quantity $\sum_{l=1}^k\{\Ex[Y \basW_l(W)] -\sum_{j=1}^k
\beta_j\Ex[\basZ_j(Z)\basW_l(W)]\}^2$ over all vectors
$(\beta_j)_{j=1}^k$. Then,~$\sol_k$ converges to the true structural
function as~$k$ tends to infinity only under an additional assumption
(e.g.\ the <<extended link condition>> introduced below) on the basis $\{\basW_j\}_{j\geq1}$.  We
are going to develop a least squares estimator $\hsol_k$ of $\sol$
based on dimension reduction and thresholding, and we show that it can
attain optimal rates of convergence in terms of a weighted risk --
provided the choice of the dimension parameter~$k$ is made in the
optimal way. It is worth to note that all the results in this article
are obtained without any additional smoothness assumption on the joint
density of $(Y,Z,W)$. In fact, such a density need not even exist.

Our main contribution is the development of a method to choose the
dimension parameter $k$ in a fully data driven way, that is, not
depending on characteristics of $\sol$, and assuming only that the
underlying conditional expectation  is <<smoothing>> in a
sense to be made precise below. The central result of the present
paper states that for this automatic choice~$\whk$, the least
squares estimator $\hsol_{\whk}$ can attain the lower bound up to a
constant, and is thus minimax-optimal. The adaptive choice of~$k$ is
made following the general model selection methodology which has been
developed in~\cite{BBM:99}. More specifically,~$\hk$ is the minimizer
of a penalized contrast. We illustrate all of our results by
considering the estimation of derivatives of the structural function
under a smoothing conditional expectation.  Typically, one
distinguishes the finitely and infinitely smoothing case.
\cite{LM:09} propose an adaptive estimator for finitely smoothing case. They derive oracle
inequalities and obtain convergence rates which differ from the
optimal ones by a logarithmic factor.  In contrast to this, we provide
a unified estimation procedure which can attain minimax-optimal rates
in either of the both cases. In other words, our estimation procedure
attains optimal rates without knowing in advance if we are in the 
finitely or infinitely smoothing case.


This article proceeds as follows. In Section~\ref{D:sec:minimax}, we
develop the minimax theory for the nonparametric instrumental
regression model \wrt the weighted risk.  We derive, as an
illustration, the optimal convergence rates for the estimation of
derivatives in the finitely and in the infinitely smoothing case.
Finally, in Section~\ref{D:sec:adapt-estim-under}, we construct the
adaptive estimator. An upper risk bound is shown and convergence rates
for the finitely and infinitely smoothing case are found to coincide
with minimax optimal ones.  The proofs and some auxiliary results are deferred to the
appendix.


\section{Minimax optimal estimation}\label{D:sec:minimax}
In this section, we develop a minimax theory for the estimation of the
structural function and its derivatives in nonparametric instrumental
regression models.

\subsection{Basic model assumptions}
It is convenient to rewrite the model equations (\ref{D:model:NP}--\ref{D:model:NP2}) in
terms of an operator between Hilbert spaces. Therefore, let us first
introduce the Hilbert spaces
\begin{gather*}
  L^2_Z=\big\{\phi:\RR^p\to \RR\;\big|\; \normV{\phi}^2_Z:=\Ex[\phi^2(Z)]<\infty\big\},\\
  L^2_W=\big\{\psi:\RR^q\to \RR\;\big|\;
  \normV{\psi}^2_W:=\Ex[\psi^2(W)]<\infty\big\},
\end{gather*}
endowed with inner products
$\skalarV{\phi,\tilde\phi}_Z=\Ex[\phi(Z)\tilde\phi(Z)]$,
$\phi,\tilde\phi\in L^2_Z$, and
$\skalarV{\psi,\tilde\psi}_W=\Ex[\psi(W)\tilde\psi(W)]$,
$\psi,\tilde\psi\in L^2_W$, respectively.  The conditional
expectation of $Z$ given $W$ defines a linear operator
$\Op\phi:=\Ex[\phi(Z)|W]$, $\phi\in L^2_Z$, which maps $L^2_Z$
to~$L^2_W$.  Taking the conditional expectation \wrt the instrument
$W$ on both sides in equation~\eqref{D:model:NP} yields
\begin{equation*}\label{D:bm:mequ}
  \gf:=\Ex[Y|W] = \Ex[\sol(Z)|W]=:\Op \varphi,
\end{equation*}
where the function $\gf$ belongs to $L^2_W$.  The estimation of the
structural function~$\sol$ is thus linked to the inversion of the
conditional operator~$\Op$. Moreover, we suppose throughout this work
that the operator $\Op$ is compact, which is the case under fairly
mild assumptions. For example, if the triple~$(Y,Z,W)$ has a joint
density, it is sufficient to demand that it be square integrable -- or
continuous, if its support is compact -- in order for $\Op$ to be
compact \citep[c.f.][]{CFR:07}.  Consequently, unlike in a
multivariate linear instrumental regression model, a continuous
generalized inverse of $\Op$ does not exist as long as the range of
the operator $\Op$ is an infinite dimensional subspace of $L^2_W$.
This corresponds to the setup of statistical ill-posed inverse
problems with unknown operator.  For a detailed discussion in the
context of inverse problems see Chapter 2.1 in \cite{EHN:96}, while in
the special case of a nonparametric instrumental regression we refer
to \cite{CFR:07}. In what follows, we always
assume that the joint distribution of $(Y,Z,W)$ is such that
$g=\E[Y|W]$ lies in the range of $\Op$ and that~$\Op$ is
injective.

\subsection{Complexity of the problem: a lower bound}
In this section we show that the obtainable accuracy of any estimator
of the structural function $\sol$ is essentially determined by
additional regularity conditions imposed on $\sol$ and the conditional
expectation operator $\Op$. In the present paper, these conditions are
characterized through different weighted norms in $L^2_Z$ with respect
to a pre-specified orthonormal basis $\{\basZ_j\}_{j\geq1}$ of
$L^2_Z$.  We formalize these conditions as follows.

\subsubsection*{Minimal regularity conditions} Given a strictly positive
sequence of weights $\beta:=(\beta_j)_{j\geqslant1}$, we denote by $\norm_\beta$
the weighted norm given by
\begin{equation*}\label{eq:weighted_norm}
  \normV{f}_\beta:=  \sum_{j=1}^\infty \beta_j |\skalarV{f,\basZ_j}_Z|^2, \qquad \forall f\in L^2_Z.
\end{equation*}
We shall measure the accuracy of any estimator $\hsol$ of the unknown
structural function in terms of a weighted risk, that is
$\Ex\normV{\hsol-\sol}_\hw^2$, for a pre-specified sequence of weights
$\hw:=(\hw_j)_{j\geqslant1}$. This general approach allows as to
consider not only the estimation of the structural function itself but
also of its derivatives, as we will explain in
Section~\ref{D:sec:minimax:ill} (illustrations) below.  Moreover,
given a sequence of weights $\solw:=(\solw_j)_{j\geqslant 1}$ we
suppose, here and subsequently, that for some constant $\solr>0$ the
structural function~$\sol$ belongs to the ellipsoid
\begin{equation*}\label{D:def:ell}
  \cF_{\solw}^\solr := \Bigl\{f\in L^2_Z\;\big|\;\normV{f}_{\solw}^2\leq \solr\Bigr\},
\end{equation*}
which captures all the prior information (such as smoothness) about
the unknown structural function $\sol$.  Furthermore, as usual in the
context of ill-posed inverse problems, we specify the mapping
properties of the conditional expectation operator~$\Op$.
More precisely, we are going to impose restrictions on the decay of
the sequence $(\normV{\Op \basZ_j}_W)_{j\geqslant1}$. Denote by $\cT$
the set of all injective compact operator mapping $L^2_Z$ to~$L^2_W$.
Given a strictly positive sequence of weights
$\Opw:=(\Opw_j)_{j\geqslant 1}$ and a constant $\Opd\geqslant 1$, we
define the subset~$\cT_{\tw}^\td$ of~$\cT$ by
\begin{equation}\label{D:bm:link}
  \cT_{\tw}^{\td}:=\Bigl\{ T\in\cT \;\big|\;   \normV{f}_{\Opw}^2/\Opd\leqslant \normV{\Op f}^2_W\leqslant {\Opd}\, \normV{f}_{\Opw}^2,\quad \forall f \in L^2_Z\Bigr\}.
\end{equation}
Notice that for all $T\in\cT_{\tw}^\td$ it follows that $\Opd^{-1}\leq
\normV{\Op\basZ_j}_W^2/\Opw_j\leq\Opd$.  Furthermore, let us denote by
$T^*:L^2_W\to L^2_Z$ the adjoint of $T$ which satisfies
$T^*\psi=\Ex[\psi(W)|Z]$ for all $\psi\in L^2_W$. 
One can show that the
sequence~$\Opw$ specifies in particular the decay of the eigenvalues of $T^*T$.  All
results of this work are derived under regularity conditions on the
structural function $\sol$ and the conditional expectation operator
$\Op$ described by the sequences $\bw$ and $\tw$, respectively.
However, below we provide illustrations of these conditions by
assuming a \flqq regular decay\frqq\ of these sequences.  The next
assumption summarizes our minimal regularity conditions on these
sequences.
\begin{assumption}\label{D:ass:minreg} Let $\solw:=(\solw_j)_{j\in \N}$,
  $\hw:=(\hw_j)_{j\in \N}$ and $\Opw:=(\Opw_j)_{j\in \N}$ be strictly
  positive sequences of weights with $\solw_0=\hw_0=\Opw_0 =1$
  such that $(\hw/\solw)$, $(\Opw/\hw)$, and $\Opw$ are
  non-increasing, respectively and such that $\zeta
  :=\sup_{k\in\N}k^3/\solw_k<\infty$, implying in particular
  $\Gamma:=\sum_{j\in\N}\solw_j^{-1}<\infty$.
\end{assumption}
It is worth noting that the monotonicity
assumption on $(\hw/\solw)$ only ensures that $\normV{\sol}_\hw$ is
finite for all~$\sol\in\Fgr$, and hence the weighted risk is a
well-defined measure of accuracy for estimators of
$\sol$. Heuristically, this reflects the fact that we cannot estimate
the $(s+1)$-th derivative if the structural function has only $s$
derivatives.

\subsubsection*{The lower bound} The next assertion provides a lower
bound for the risk \wrt the weighted norm. Thus, we extend the result
of \cite{CR:10}, who show a lower bound for the mean integrated
squared error.
\begin{theorem}\label{D:res:lower} Suppose that the
  iid.\ $(Y,Z,W)$-sample of size $n$ obeys the model
  (\ref{D:model:NP}--\ref{D:model:NP2}), that the distribution of the error term $U$ belongs
  to the class
 \[\cU_\sigma:=\{ P_U \;|\; \Ex [U|W] =0 \mbox{ and }\Ex [U^4|W] \leq
  \sigma^4\}\] with $\sigma>0$ and that $\sup_{j\geq1} \Ex [e_j^4(Z)
  |W]\leq \eta$, $\eta\geq1$. Consider sequences $\solw$, $\hw$ and
  $\Opw$ satisfying Assumption~\ref{D:ass:minreg} such that the conditional
  expectation operator~$T$ associated to $(Z,W)$ belongs to
  $\cT_{\tw}^\td$, $\td\geq 1$.  Define for all $n\geq 1$
  \begin{multline}\label{D:def:kn}
    \kstar:=\kstar(\solw,\Opw,\hw):=\argmin\limits_{k\in\N}\Bigl\{\max\Bigl(\frac{\hw_k}{\solw_k},
    \sum_{j=1}^{k}\frac{\hw_j}{n\Opw_j}\Bigr)\Bigr\}
    \mbox{ and }\\
    \dstar:=\dstar(\solw,\Opw,\hw):=\max\Bigl(
    \frac{\hw_{\kstar}}{\solw_{\kstar}},
    \sum_{j=1}^{\kstar}\frac{\hw_j}{n\Opw_j}\Bigr).
  \end{multline}
  If in addition
  $\kappa:=\inf_{n\geq1}\{(\dstar)^{-1}\min(\hw_{\kstar}\solw_{\kstar}^{-1},
  \sum_{l=1}^{\kstar}{\hw_l}{(n\Opw_l})^{-1}) \}>0$ and $ \sigma^4\geq
  8(3+2\solr^2 \Gamma^2)$, then for all $n\geq 1$ and for any
  estimator $\tsol$ of $\sol$, we have
  \begin{equation*} \sup_{P_U\in \cU_\sigma} \sup_{\sol
      \in\cF_\solw^{\solr}} \Ex\normV{\tsol-\sol}^2_\hw \geqslant
    \frac{\kappa}{4}\,\min\bigg(\solr, \frac{1}{2\Opd}\bigg)\,\dstar.
  \end{equation*}
\end{theorem}

\begin{remark}\label{D:rem:lower} 
  The proof of the last assertion is based on Assouad's cube technique
  \citep[c.f.][]{KorostelevTsybakov1993,Tsy:04}, which consists in
  constructing $2^{\kstar}$ candidates of structural functions which
  have the largest possible $\norm_\hw$-distance but are still
  statistically non distinguishable.  In the last theorem, the
  additional moment condition $\sup_{j\geq1} \Ex [e_j^4(Z) |W]\leq
  \eta$ is obviously satisfied if the basis functions $\{\basZ_j\}$
  are uniformly bounded (e.g.\ the trigonometric basis considered in
  Section~\ref{D:sec:minimax:ill}). However, if $V$ denotes a Gaussian
  random variable with mean zero and variance one, which is moreover
  independent of $(Z,W)$, then the additional condition $\sigma^4\geq
  8(1+2\solr^2 \Gamma^2\eta)$ ensures that for all structural
  functions $\sol \in\cF_\bw^{\br}$, the distribution of the error
  term $U:= V -\sol(Z) + [\Op\sol](W)$ belongs to $\cU_\sigma$. This
  specific case is only needed to simplify the calculation of the
  distance between distributions corresponding to different structural
  functions. A similar assumption has been used by \cite{CR:10}. 

  On the other hand, below we derive an upper bound assuming that the
  distribution of error term $U$ belongs to $\cU_\sigma$ and that the
  joint distribution of $(Z,W)$ satisfies additional moment
  conditions. In this situation, Theorem~\ref{D:res:lower} provides a
  lower bound for any estimator as long as $\sigma$ is sufficiently
  large.  Note further that this lower bound tends only to zero if
  $\hw/\solw$ is a vanishing sequence. In other words, in case
  $\solw\equiv
  1$, 
  uniform consistency over all $\sol$ with $\normV{\sol}_Z^2\leq
  \solr$ can only be achieved with respect to a weighted norm weaker
  than the $L_Z^2$-norm, that is, if $\hw$ is a sequence tending to
  zero. Finally, it is important to note that the regularity
  conditions imposed on the structural function~$\sol$ and the
  conditional expectation operator~$\Op$ involve only the basis
  $\{\basZ_j\}_{j\geq1}$ in $L^2_Z$. Therefore, the lower bound
  derived in Theorem~\ref{D:res:lower} does not capture the influence
  of the basis $\{\basW_l\}_{l\geq1}$ in~$L^2_W$ used to construct the
  estimator. In other words, the proposed estimator of $\sol$ can only
  attain this lower bound if $\{\basW_l\}_{l\geq1}$ is appropriately
  chosen.\qed
\end{remark}

\subsection{Minimax-optimal Estimation by dimension reduction and thresholding}
In addition to the basis $\{\basZ_j\}_{j\geq1}$ of $L^2_Z$ considered
in the last section, we introduce now a basis $\{\basW_l\}_{l\geq1}$
in $L^2_W$.  In this section we derive the asymptotic properties of
the least squares estimator under minimal assumptions on these two
bases. More precisely, we suppose that the structural function $\sol$
belongs to some ellipsoid $\cF_\bw^\br$ and that the conditional
expectation satisfies a link condition, i.e., $T\in\cTdw
$. Furthermore, we introduce an additional condition linked to the
basis $\{\basW_l\}_{l\geq1}$. Then we show that the proposed estimator
attains the lower bound derived in the last section. All these results
are illustrated under classical smoothness assumptions at the end of
this section.

\subsubsection*{Matrix and operator notations} Given $k\geq 1$, let
$\cE_k$ and $\cF_k$ denote the subspace of $L^2_Z$ and $L^2_W$ spanned
by the functions $\{e_j\}_{j=1}^k$ and $\{f_l\}_{l=1}^k$,
respectively. $E_k$ and $E_k^\perp$ (resp. $F_k$ and $F_k^\perp$)
denote the orthogonal projection mappings on $\cE_k$ (resp. $\cF_k$)
and its orthogonal complement~$\cE_k^\perp$ (resp. $\cF_k^\perp$),
respectively. Given a matrix $K$, its inverse is denoted by $K^{-1}$ and its
transposed matrix by $K^t$. Let $[\phi]$, $[\psi]$ and $[K]$
denote the (infinite) vector and matrix of the function $\phi\in
L^2_Z$, $\psi\in L^2_W$ and the operator $K:L^2_Z\to L^2_W$ with the
entries $[\phi]_{j}=\skalarV{\phi,e_j}$,
$[\psi]_{l}=\skalarV{\psi,f_l}$ and $[K]_{lj}=\skalarV{Ke_j,f_l}$,
respectively. The upper $k$-sub-vector and $(k\times k)$-sub-matrix of
$[\phi]$, $[\psi]$ and $[K]$ are denoted by $[\phi]_{\uk}$,
$[\psi]_{\uk}$ and $[K]_{\uk}$, respectively. Note that
$[K^\ast]_{\uk}=[K]_{\uk}^t$. The diagonal matrix with entries~$v$ is
denoted by $\diag(v)$ and the identity matrix is denoted by
$I$.  Clearly, $[E_k \phi]_{\uk} =[\phi]_{\uk}$ and if we restrict
$F_k K E_k$ to an operator from $\cE_k$ into $\cF_k$, then it has the
matrix $[K]_{\uk}$. Moreover, if $v\in\R^k$ then $\normV{v}$ denotes
the Euclidean norm of $v$, and given a $(k\times k)$-matrix~$M$, let
$\normV{M}:=\sup_{\normV{v}\leq 1}\normV{Mv}$ denote its spectral-norm
and $\tr(M)$ its trace.

Consider the conditional expectation operator $\Op$ associated to the
regressor~$Z$ and the instrument $W$. If $[\basZ(Z)]_{\uk}$ and $[\basW(W)]_{\uk}$
denote the $k$-dimensional random vectors with entries $\basZ_j(Z)$ and
$\basW_j(W)$ respectively, then $[\Op]_{\uk}=\Ex
[\basW(W)]_{\uk}[\basZ(Z)]_{\uk}^t$ which we 
assume to be non singular for all $k\geq 1$ (or, at least for
sufficiently large
 $k$), such that $[\Op]_{\uk}^{-1}$ always exists. Note that it is
a nontrivial problem to determine in under what precise conditions
such an assumption holds (see e.g.\ \cite{EK:01}
and references therein).

\subsubsection*{Definition of the estimator}Let
$(Y_1,Z_1,W_1),\dotsc,(Y_n,Z_n,W_n)$ be an iid.\ sample of $(Y,Z,W)$.
Since $[\Op]_{\uk}=\Ex [\basW(W)]_{\uk}[\basZ(Z)]_{\uk}^t$ and
$[\gf]_{\uk}=\Ex{Y[\basW(W)]_{\uk}}$ can be written as expectations,
we can construct estimators by using their empirical counterparts,
that is,
\begin{equation*}
\hfou{\Op}_{\uk}:=
(1/n)\sum_{i=1}^n [\basW(W_i)]_{\uk}[\basZ(Z_i)]_{\uk}^t \quad\mbox{
  and }\quad\hfou{\gf}_{\uk}:= (1/n)\sum_{i=1}^n Y_i[\basW(W_i)]_{\uk}.
\end{equation*}
Then the estimator of the structural function $\sol$ is defined by
\begin{equation}\label{D:gen:def:est}
  \hsol_k:= \sum_{j=1}^k[\hsol_k]_{j}\basZ_j\hfill\mbox{ with }\hfill[\hsol_k]_{\uk}:= \left\{\begin{array}{ll} 
      \hfou{\Op}_{\uk}^{-1} [\hgf]_{\uk}, &\begin{array}{l}\mbox{if $\hfou{\Op}_{\uk}$ is nonsingular}\\\mbox{and }\normV{\hfou{\Op}^{-1}_{\uk}}\leq \sqrt n,\end{array}\\[3ex]
      0,&\begin{array}{l}\mbox{otherwise},\end{array}\end{array}\right.
\end{equation}
where the dimension parameter $k=k(n)$ has to tend to infinity as the
sample size~$n$ increases. This estimator $\hsol_k$ takes its
inspiration from the linear Galerkin approach
\citep[c.f.][]{EK:01,HR:08}.

\subsubsection*{Extended link condition} Consistency of this estimator
is only possible if the least squares solution
$\sol_k=\sum_{j=1}^k[\sol_k]_{j}\basZ_j$ with
$[\sol_k]_{\uk}=[\Op]_{\uk}^{-1}[\gf]_{\uk}$ converges to $\sol$ as
$k\to\infty$, which is not true in general. However, the condition
$\sup_{k\in\N}\normV{[\Op]^{-1}_{\uk} \fou{\Op
    E_k^\perp}_{\uk}}<\infty$ is known to be sufficient to ensure
convergence of $\sol_k$. Notice that this condition involves also the
basis $\{\basW_l\}_{l\geq1}$ in $L^2_W$. In what follows, we introduce
an alternative but stronger condition to guarantee the convergence,
which extends the link condition \eqref{D:bm:link}, that is,
$T\in\cTdw$. We denote by $\cTdDw$ for some $\OpD\geq \Opd$ the subset
of $\cTdw $ given by
\begin{equation}\label{D:bm:link:gen}
  \cTdDw:=\Bigl\{ T\in \cTdw \;\big|\; \sup_{k\in\N}\normV{[\diag(\Opw)]^{1/2}_{\uk}[T]^{-1}_{\uk}}^2\leq \OpD\Bigr\}.
\end{equation}
\begin{remark}\label{D:rem:ext:link} The link condition \eqref{D:bm:link} implies the extended
  link condition \eqref{D:bm:link:gen} for a suitable $\OpD>0$ if
  $\{e_j\}$ and $\{f_j\}$ are the singular functions of $\Op$ and if
  $[\Op]$ is only a small perturbation of $\diag(\Opw^{1/2})$, or if
  $\Op$ is strictly positive (for a detailed discussion we refer to
  \cite{EK:01} and \cite{CardotJohannes2008}).
  We underline that once both bases $\{\basZ_j\}_{j\geq1}$ and
  $\{\basW_l\}_{l\geq1}$ are specified, the extended link condition
  \eqref{D:bm:link:gen} restricts the class of joint distributions of
  $(Z,W)$ to those for which the least squares solution~$\sol_k$ is
  $L^2$-consistent. Moreover, we show below that under the extended
  link condition the least squares estimator of $\sol$ given in~\eqref{D:gen:def:est} can
  attain minimax-optimal rates of convergence. In this sense, given a
  joint distribution of $(Z,W)$, a basis $\{\basW_l\}_{l\geq1}$
  satisfying the extended link condition can be interpreted as a set
  of optimal instruments.  Furthermore, for each pre-specified basis
  $\{\basZ_j\}_{j\geq1}$, we can theoretically construct a basis
  $\{\basW_l\}_{l\geq1}$ of optimal instruments such that the extended link condition is not
  a stronger restriction than the link condition~\eqref{D:bm:link} (see
  \cite{BJ:09} for more details).
  \hfill$\square$\end{remark}

\subsubsection*{The upper bound}
The following theorem provides an upper bound under the extended link
condition~\eqref{D:bm:link:gen} and an additional moment condition on
the bases or, more precisely, on the random vectors $[\basZ(Z)]$ and
$[\basW(W)]$. We begin this section by formalizing this additional
condition.
\begin{assumption}\label{D:ass:A2}
  There exists $\eta\geqslant 1$ such that the joint distribution of
  $(Z,W)$ satisfies
  \begin{itemize}
  \item[(i)] $\sup_{j\in \N}\Ex [e_j^2(Z)|W]\leqslant \eta^2$ and
    $\sup_{l\in \N}\Ex [f_l^4(W)]\leqslant \eta^4$;
  \item[(ii)] $\sup_{j,l\in\N} \var (e_j(Z)f_l(W))\leqslant \eta^2$
  \item[(iii)] $\sup_{j, l\in\N} \Ex| e_j(Z)f_l(W)- \Ex
    [e_j(Z)f_l(W)]|^k\leqslant \eta^{k-2} k! \var(e_j(Z)f_l(W))$, for all
    $k\geq 3$.
  \end{itemize}
\end{assumption}

\noindent This assumption restricts the set of
possible joint distribution of $(Z,W)$.  More precisely, it supposes
that the random variables $e_j(Z)f_l(W)- \Ex [e_j(Z)f_l(W)]$ satisfy
Cramer's condition uniformly, which is known to be sufficient to
obtain an exponential bound for their large deviations
\citep[c.f.][]{Bosq1998}.  It is however noticeable that the
assumption is satisfied for any joint distribution and for
sufficiently large $\eta$ if the bases $\{e_j\}_{j\geq1}$ and
$\{f_l\}_{l\geq1}$ are uniformly bounded.

\begin{theorem}\label{D:res:upper:A3}
  Suppose that the iid.\ $(Y,Z,W)$-sample of size $n$ obeys the model
  (\ref{D:model:NP}--\ref{D:model:NP2}) and that the joint
  distribution of $(Z,W)$ satisfies Assumption~\ref{D:ass:A2} for some
  $\eta\geq1$.  Consider sequences $\solw$, $\hw$ and $\Opw$
  satisfying Assumption~\ref{D:ass:minreg}. Let~$\kstar$,~$\dstar$,
  and $\kappa$ be as given in Theorem~\ref{D:res:lower} and suppose
  that
  \begin{equation}
    \label{eq:thm:upp}
    (\kstar)^2\;\max\left\{|\log \dstar|, (\log\kstar)\right\} =
    o(\Opw_{\kstar}) \hspace{3em} \text{as } n\to\infty.
  \end{equation}
Then, we have for all  $n\in\N$ that
\[ \sup_{\Op\in \cTdDw}   \sup_{P_U\in \cU_\sigma} \sup_{\sol \in\cF_\solw^{\solr}}
    \Ex\normV{\hsol_{\kstar}-\sol}^2_\hw \leq C\; \dstar\]
for a constant $C>0$ depending only on the classes
$\cTdDw,\cF_\solw^{\solr}$, and the constants $\sigma$ and $\eta$. 
\end{theorem}

\begin{remark}\label{D:rem:upper:A3} From Theorems~\ref{D:res:lower} and~\ref{D:res:upper:A3}
  it follows that under Assumption~\ref{D:ass:A2}, the estimator $\hsol_{\kstar}$ attains the optimal rate
  $\dstar$ for all sequences $\bw$, $\hw$ and $\tw$ satisfying the
  minimal regularity conditions from Assumption~\ref{D:ass:minreg}. Let
  us briefly discuss the role of the sequences $\bw$, $\hw$
  and $\tw$.  Theorem~\ref{D:res:lower} and~\ref{D:res:upper:A3} show that
  the faster the sequence $\tw$ decreases, the slower the obtainable
  optimal rate of convergence becomes. On the other hand, a faster
  increase of $\solw$ or decrease of $\hw$ leads to a faster optimal
  rate. In other words, as expected, a structural function satisfying
  a stronger regularity condition can be estimated faster, and
  measuring the accuracy with respect to a weaker norm leads to faster
  rates, too.  \hfill$\square$\end{remark}

\subsubsection*{Illustration: estimation of derivatives}\label{D:sec:minimax:ill}
To illustrate the previous results, we will describe in this section
the prior information about the unknown structural function $\sol$ by
its degree of smoothness. In order to simplify the presentation, we
follow \cite{HallHorowitz2005} and suppose that the marginal
distribution of the scalar regressor $Z$ and the scalar instrument~$W$
are uniformly distributed on the interval $[0,1]$. It is worth noting
that all the results below can be extended to the multivariate case in
a straightforward way. In the univariate case, it follows that both
Hilbert spaces~$L^2_Z$ and~$L^2_W$ are isomorphic to $L^2[0,1]$,
endowed with the usual norm $\norm$ and inner product $\skalar$.

In the last sections, we have seen that the choice
of the basis $\{e_j\}_{j\geq1}$ is directly linked to the a priori
assumptions we are willing to impose on the structural function.  In
case of classical smoothness assumptions, it is natural to consider
the Sobolev space of periodic functions. Therefore, we introduce the
trigonometric basis
\begin{equation*}
  \basTr_{1}:\equiv1, \;\basTr_{2j}(s):=\sqrt{2}\cos(2\pi j s),\;
  \basTr_{2j+1}(s):=\sqrt{2}\sin(2\pi j s),s\in[0,1],\; j\in\N.
\end{equation*}
and choose $\{e_j=\psi_j\}$. It is well-known that for a weight
sequence $\solw$ with $\solw_1=1$ and $\solw_j=j^{2p}$ for $j\geq2$,
the ellipsoid $\Fgr$ is a subset of the Sobolev space of $p$-times
differentiable periodic functions.
In the rest of this section we will suppose that the prior information
about the unknown structural function $\sol$ is characterized by such
a Sobolev ellipsoid, i.e.\ that $\sol$ is $p\geq 0$ times differentiable. In
this illustration, we consider the estimation of derivatives of the
structural function $\sol$. We therefore recall
that, up to a constant, for any function $h\in\Fgr$ the
weighted norm $\normV{h}_\hw$ with $ \hw_0=1\mbox{ and }
\hw_{j}=j^{2s},$ $j\geq 2$, equals the $L^2$-norm of the $s$-th weak
derivative $h^{(s)}$ for each integer $0\leq s\leq p$. By virtue of
this relation, the results in the previous section imply also a lower
as well as an upper bound of the $L^2$-risk for the estimation of the
$s$-th weak derivative of $\sol$. Finally, we restrict our attention
to conditional expectation operator $\Op\in\cTdw$ with either
\begin{itemize}
\item[{\bf[p-$\pmb{\Opw}$]}] a polynomially decreasing  sequence
  $\Opw$, i.e., $\Opw_0=1$ and $\Opw_j =
  j^{-2a}$, $j\geq2$, for some $a>0$, or
\item[{\bf[e-$\pmb{\Opw}$]}] an exponentially decreasing  sequence
  $\Opw$, i.e., $\Opw_0=1$ and $\Opw_j =
  \exp(-j^{2a})$, $j\geq2$, for some $a>0$.
\end{itemize}
It is easily seen that the minimal regularity conditions given in
Assumption~\ref{D:ass:minreg} are satisfied if $p>1/2$.  Roughly
speaking, this means that the structural function is at least
continuous.  The lower bound presented in the next assertion follows
now directly from Theorem~\ref{D:res:lower}. Note that the additional
condition, $\sup_{j\geq1} \Ex [e_j^4(Z) |W]\leq \eta$, $\eta\geq8$, is
satisfied since the trigonometric basis is bounded uniformly by two.
Before stating the results, let us introduce some
asymptotic notation: We write\label{C:sim:sim:sim} $a_n\lesssim b_n$
when there is a $C\in\RR_+$ such that $a_n\leqslant C\, b_n$ for all sufficiently large $n\in\N$ and
$a_n\sim b_n$ when $a_n\lesssim b_n$ and $b_n\lesssim a_n$
simultaneously.

\begin{proposition}\label{D:coro:ex:lower} Suppose an   iid.\ sample of size $n$ from the model
  (\ref{D:model:NP}--\ref{D:model:NP2}).  If $\solw_j=j^{2p}$ with $p>1/2$,
  then we have for any
  estimator $ \tsol^{(s)}$ of $\sol^{(s)}$, $0\leq s<p$, \\[-.5em]
  \begin{itemize}\item[{\bf[p-$\pmb{\Opw}$]}] 
     \hspace*{5ex}$ \sup_{P_U\in \cU_\sigma}\sup_{\sol\in\Fgr} \left\{ \Ex\normV{\tsol^{(s)}-\sol^{(s)}}^2\right\}\gtrsim
    n^{-2(p-s)/(2p+2a+1)},$
  \item[{\bf[e-$\pmb{\Opw}$]}] 
  \hspace*{5ex}$  \sup_{P_U\in \cU_\sigma}\sup_{\sol\in\Fgr} \left\{ \Ex\normV{\tsol^{(s)}-\sol^{(s)}}^2\right\}\gtrsim
    (\log n)^{-(p-s)/a}$.
  \end{itemize}\end{proposition}
\vspace{1em}

\noindent In this section, the basis of $L^2_W$ is given by the
trigonometric basis $\{\basW_l=\basTr_l\}_{l\geq1}$.  The moment
conditions formalized in Assumption~\ref{D:ass:A2} are thus
automatically fulfilled since the bases $\{\basZ_j\}_{j\geq1}$ and
$\{\basW_l\}_{l\geq1}$ are both uniformly bounded. We suppose that the
associated conditional expectation operator $\Op$ satisfies the
extended link condition \eqref{D:bm:link:gen}, that is,
$\Op\in\cTdDw$. By this means, we restrict the set of possible joint
distributions of $(Z,W)$ to those having the trigonometric basis as
optimal instruments. As an estimator of $\sol^{(s)}$, we shall
consider the $s$-th weak derivative of the estimator $\hsol_k$ defined
in \eqref{D:gen:def:est}. Recall that for each integer $0\leq s\leq
p$, the $s$-th weak derivative of the estimator $\hsol_k$ is
\begin{equation*}
  \hsol^{(s)}_k(t)=\sum_{j\in\Z}(2i \pi j)^s\int_0^1 \hsol_k(u) \exp(-2i\pi ju)du \exp(-2i\pi j t).
\end{equation*}
Applying Theorem~\ref{D:res:upper:A3}, the rates of the lower bound given
in the last assertion are seen to coincide, up to a constant, with an
upper bound of the $L^2$-risk of the estimator $\hsol^{(s)}_k$, which
is the statement of the next proposition. This proves that these rates
are optimal and the estimator $\hsol^{(s)}_k$ is minimax optimal in
both cases.

\begin{proposition}\label{D:coro:ex:upper} Suppose that the
  iid. $(Y,Z,W)$-sample of size $n$ obeys the model
  (\ref{D:model:NP}--\ref{D:model:NP2}).  Let $\solw_j=j^{2p}$ for
  $p\geq 3/2$.
  For $0\leq s<p$ consider the estimator  $\hsol_{\kstar}$ given in \eqref{D:gen:def:est}.\\[-4ex]
  \begin{itemize}\item[{\bf[p-$\pmb{\Opw}$]}] In the polynomial decreasing case with  $\kstar\sim n^{1/(2p+2a+1)}$,\\[1ex]
    \hspace*{5ex}$\sup_{P_U\in \cU_\sigma}\sup_{\sol\in\Fgr} \left\{
      \Ex\normV{\hsol^{(s)}_{\kstar}-\sol^{(s)}}^2\right\}\lesssim
    n^{-2(p-s)/(2p+2a+1)}$.
  \item[{\bf[e-$\pmb{\Opw}$]}] In the exponentially decreasing case  with $\kstar\sim (\log n)^{1/(2a)}$,\\[1ex]
    \hspace*{5ex}$\sup_{P_U\in \cU_\sigma}\sup_{\sol\in\Fgr } \left\{
      \Ex\normV{\hsol^{(s)}_{\kstar}-\sol^{(s)}}^2\right\}\lesssim
    (\log n)^{-(p-s)/a}$.
  \end{itemize}
\end{proposition}

 \begin{remark}\label{D:rem:upper:sob:1}
   We emphasize the interesting role of the parameters $p$ and $a$
   characterizing the regularity conditions imposed on $\sol$ and
   $\Op$ respectively: As we see from Propositions~\ref{D:coro:ex:lower}
   and~\ref{D:coro:ex:upper}, if the value of $a$ increases, the
   obtainable optimal rate of convergence decreases. Therefore, the
   parameter $a$ is often called {\it degree of ill-posedness}
   \citep[c.f.][]{Natterer84}.  On the other hand, an increase of the
   quantity~$p$ leads to a faster optimal rate. In other words, as
   expected, a smoother structural function can be estimated faster.
   Finally, as opposed to the polynomial case, in the exponential case
   the smoothing parameter $\kstar$ does not depend on the value of~$p$. It follows that the proposed estimator is automatically
   adaptive, i.e.\ it does not depend on an a-priori knowledge of the
   degree of smoothness of the structural function $\sol$. However,
   the choice of the smoothing parameter does depend on the properties
   of $\Op$, more precisely, the value of
   $a$.\hfill$\square$\end{remark}


\section[Adaptive estimation]{Adaptive estimation under smoothness assumptions}
\label{D:sec:adapt-estim-under}
In this section, our objective is to construct a fully adaptive
estimator of the structural function $\sol$.  Adaptation means that in
spite of the conditional expectation operator $\Op$ being unknown, the
estimator should attain the optimal rate of convergence over the
ellipsoid $\cF_\solw^\solr$ for a wide range of different weight
sequences $\solw$. However, we will suppose that the operator $\Op$ is
diagonal with respect to the trigonometric basis $\{\basTr_j\}$. In
this situation, for example, an operator with polynomially
decreasing~$\Opw$ having a degree of ill-posedness $a$ behaves like
$a$-times integrating, and hence it is also called \textit{finitely
  smoothing}. On the other hand, when the sequence~$\Opw$ is
exponentially decreasing with degree of ill-posedness $a$, the
operator behaves like integrating infinitely many times, and hence it
is also called \textit{infinitely smoothing}. Thus, this additional
condition imposes in fact a smoothing condition on the unknown
conditional expectation operator $\Op$. Even though we assume that the
operator is smoothing, we do not impose any a-priori knowledge about
the specific decay of $\Opw$.  Our starting point is the estimator
given in \eqref{D:gen:def:est}, which in this situation takes the form
\begin{equation}
  \hsol_k =
  \sum_{j=1}^k\frac{\hfou{\gf}_j}{\hfou{\Op}_{jj}}\I{[\inf_{1\leq j \leq k}\hfou{\Op}_{jj}^2
    \geq 1/n]}\;\basTr_j,\label{D:eq:17}
\end{equation}
with $\hfou{\gf}_j$ and $\hfou{\Op}_{jj}$ defined in
\eqref{D:gen:def:est}. In the last section, we have shown that this
estimator is minimax-optimal provided the dimension parameter $k$ is
chosen in the optimal way.  In what follows, the dimension parameter
$k$ is chosen using a model selection approach via penalization. This
choice will only involve the data and none of the sequences $\solw$
and $\Opw$ describing the underlying smoothness. First, we introduce
some sequences which are used below.

\enlargethispage{2em}
\begin{definition}\label{D:def:known}
  \ii\
  \begin{enumerate}\item For all $k\geq 1$, define
    $\Delta_k:=\max_{1\leq j\leq k}\hw_j/\Opw_j$, $\tau_k:=\max_{1\leq
      j\leq k}(\hw_j)_{\vee1}/\Opw_j$ with $(q)_{\vee1}:=\max(q,1)$
    and
    \begin{equation*}\label{D:deltamref} \delta_k := {k}
      \Delta_k\frac{\log (\tau_k\vee (k+2))}{\log(k+2)}.
    \end{equation*}
    Let further~$\Sigma$ be a non-decreasing function such that for
    all $C>0$
    \begin{equation}\label{D:ass:sum}
      \sum_{k\geq
        1}C\,{\tau_k}\exp\Big(-\frac{k\log(\tau_k \vee (k+2))}{6C\log(k+2)}\Big)\leq
      \Sigma(C)<\infty
    \end{equation}
    and $ \sup_{n\in\N}\,
    \exp\big(-K_2\;C^{-1}\;\,n^{1/6}+\frac{5}{3}\log n\big)\leq
    \Sigma(C)$ with the constant $K_2 = (\sqrt{2} -1 )/(21\sqrt{2})$.
  \item Define a sequence $\Cy$ as follows,
    \begin{multline*}
      \Cy_n := \Cy_n(\Opw,\Opd) := \max\, \bigg\{1\leq N \leq n
      \;\bigg|\; n^7\,\exp\Big(-\frac{n\,\Opw_N}{288\Opd}\Big) \leq
      \Big(\frac{2016\,\Opd}{\Opw_1}\Big)^7\\\text{\rm and}\quad \delta_N
      / n \leq 1 \bigg\}.
    \end{multline*}
  \end{enumerate}\arab
\end{definition}
\noindent It is easy to see that there exists always a function
$\Sigma$ satisfying condition \eqref{D:ass:sum}. Consider the estimator
$\hsol_{\wtk}$ defined by choosing the dimension parameter $\wtk$ such
that
\begin{equation*}
  \wtk := \argmin_{1\leq k\leq \Cy_n}\left\{-\normV{\hsol_{k}}^2_\hw +
    c\,\frac{\delta_k}{n}  \right\}
\end{equation*}
for some constant $c>0$.  
However, the estimator $\hsol_{\wtk}$ is only partially adaptive,
since the dimension parameter is chosen using a criterion function
that involves the sequences $\Cy$ and $\delta$ which depend on $\Opw$
and $\Opd$. We solve this problem by defining empirical versions of
these sequences. The fully adaptive estimator is then defined
analogously to the one above, but uses the estimated rather than the
original sequences.

\begin{definition}\label{D:def:unknown} Let $\hdelta:=(\hdelta_k)_{k\geq1}$, $\hCy:=(\hCy_n)_{n\geq1}$, be as follows.
  \ii\begin{enumerate}\item Given $\hDelta_k := \max_{1\leq j\leq
      k}\hw_j \hfou{\Op}_{jj}^{-2}\I{[\inf_{1\leq j\leq
        k}\hfou{\Op}_{jj}^2\geq 1/n]}$ and\\ $\htau_k :=\max_{0\leq
      j\leq k} (\hw_j)_{\vee1} \hfou{\Op}_{jj}^{-2}\I{[\inf_{1\leq
        j\leq k}\hfou{\Op}_{jj}^2\geq 1/n]}$ let
    \[\hdelta_k := k \hDelta_k \,\frac{\log(\htau_k \vee
      (k+2))}{\log(k+2)}.\]
  \item Given $\Cy_n^u := \argmax_{1\leq N\leq n}\big\{\max_{1\leq
      j\leq N} {\hw_j}/{n}\leq 1\big\}$, let
    \[\hCy_n := \quad\argmin_{1\leq j \leq \Cy_n^u}
    \bigg\{\frac{|\hfou{\Op}_j|^2}{|j|(\hw_j)_{\vee1}} < \frac{\log
      n}{n}\bigg\}.\]

  \end{enumerate}
\end{definition}
It worth to stress that all these sequences do not involve any
a-priori knowledge about neither the target function $\sol$ nor the
operator $\Op$. Now, we choose the dimension parameter as
\begin{equation}
  \hk := \argmin_{1\leq k\leq \hCy_n} \bigg\{ -\normV{\hf_k}^2_\hw  +
  540\,\Ex[Y^2]\,\frac{\hdelta_k}{n}  \bigg\}.\label{D:eq:21}
\end{equation}
Throughout this chapter we do not address the issue that the value
$\Ex[Y^2]$ is not known in practice. Anyway, it can easily be
estimated by its empirical counterpart.  Moreover, the constant $540$,
though suitable for the theory, may probably be chosen much smaller in
practice by a simulation study (cf.\ \cite{CRT:06.2} in the context of
a deconvolution problem).

\bigskip
\noindent Our main result below requires the following Assumption.
\begin{assumption}\label{D:ass:known}\label{D:ass:unknown}\ii
  $\,$ The sequence $\Cy$ from Definition~\ref{D:def:known}~(ii)
    satisfies the conditions
    \begin{equation*}
      \max_{j\geq\Cy_n}\frac{\Opw_j}{j(\hw_j)_{\vee1}}\leq\frac{\log
        n}{4 \Opd n }
      \qquad\text{ and }\qquad \edfr^{-1}\,\min_{1\leq j\leq\Cy_n}
      \Opw_j \geq 2/n. 
    \end{equation*}
\end{assumption}

\begin{remark}\label{D:sec:adapt-estim-under-1}
  Assumption~\ref{D:ass:unknown} is satisfied for sufficiently large $n$
  by construction. Let us illustrate briefly this assumption in the
  setting of the examples introduced in
  Section~\ref{D:sec:minimax:ill}.  Recall the distinction between
  finitely and infinitely smoothing conditional expectation operators.
  The sequences from Definition~\ref{D:def:known} take the following
  forms in the two respective cases.
  \begin{itemize}
  \item[{\bf[fs]}] In the finitely smoothing case, we have
    \[ \Delta_k = k^{2a+2s},\quad \delta_k \asymp
    k^{2a+2s+1},\quad\Cy_n \asymp
    n^{1/(2a+2s+1)}.\]
  \item[{\bf[is]}] In the infinitely smoothing case, we have
    \begin{multline*}
      \Delta_k = k^{2s}\exp(k^{2a}),\quad  \delta_k \asymp k^{2a+2s+1}\exp(k^{2a})(\log k)^{-1},\\
      \Cy_n \asymp \bigg(\log\frac{ n \;\log \log n}{(\log
        n)^{(2a+2s+1)/(2a)}}\bigg)^{1/(2a)}.
    \end{multline*}
  \end{itemize}
 The sequence $N$ satisfies Assumption~\ref{D:ass:unknown} in either case.\qed
\end{remark}
We are now able to state the main result of this chapter providing an
upper risk bound for the fully adaptive estimator in the case where
the eigenfunctions of the operator $\Op^*\Op$ are known.
\begin{theorem}\label{D:thm:upper:unknown}
  Assume that we have a sample of size $n$ of $(Y,Z,W)$. Consider
  sequences $\hw,$ $\solw$, and $\Opw$ satisfying
  Assumption~\ref{D:ass:minreg} such that the conditional expectation
  operator $T$ associated to $(Z,W)$ belongs to $T\in\cTdDw$,
  $\Opd,\OpD\geq1$ and is diagonal \wrt the trigonometric basis
  $\{\psi_j\}$. Let the sequences~$\delta$ and~$\Cy$ be as in
  Definition~\ref{D:def:known} and suppose that
  Assumption~\ref{D:ass:unknown} holds.  Define further $\Cy^l_n :=
  \argmax_{1\leq j \leq \Cy_n}\big\{ \frac{\Opw_j}{j(\hw_j)_{\vee1}}
  \geq \frac{4\Opd\,\log n}{n}\big\}$.  Consider the estimator
  $\hsol_{\hk}$ defined in~\eqref{D:eq:17} with $\whk$ given
  by~\eqref{D:eq:21}.  Then for all $n\geq 1$
  \begin{multline*}
    \sup_{P_U\in \cU_\sigma} \sup_{\sol \in\cF_\solw^{\solr}} \left\{
      \Ex\normV{\hsol_{\hk}-\sol}_\hw^2\right\} \lesssim
    (2\solr\Gamma+\sigma^2+1)^4 \Opd\,\zeta_\Opd \bigg[\min_{1\leq k
      \leq \Cy_n^l}\bigg\{\max\bigg(\frac{\hw_k}{\solw_{k}},
    \frac{\delta_{k}}{n}\bigg)\bigg\} \\ + \solr \max_{j\geq 1}
    \bigg\{
    \frac{\hw_j}{\solw_j}\min\bigg(1,\frac{1}{n\Opw_j}\bigg)\bigg\}
     + \frac{1}{n}\, \,
    \bigg\{\Sigma\left(\frac{(2\solr\Gamma+\sigma^2)\zeta_\Opd+\vuz}{\vuz^2}\right)+1\bigg\}
    \bigg],
  \end{multline*}
  where $\vuz :=\Ex[\var(U|Z)]$ and $\zeta_d := (\log 3d) / \log 3$.
\end{theorem}
\noindent Compare the last assertion with the lower bound given in
Theorem~\ref{D:res:lower}. It is easily seen that if $(\hw/\Opw)$ is
non-decreasing, the second term in the upper bound of
Theorem~\ref{D:thm:upper:unknown} is always smaller than the first
one. Thus, in this situation the fully adaptive estimator attains the
lower bound up to a constant if and only if
\[R_n^\diamond := \min_{1\leq k \leq
  \Cy_n^l}\bigg\{\max\bigg(\frac{\hw_k}{\solw_{k}},
\frac{\delta_{k}}{n}\bigg)\bigg\}\] is of the same order as $R_n^\ast
= \min_{k\in\N}\Bigl\{\max\Bigl(\frac{\hw_k}{\solw_k},
\sum_{j=1}^{k}\frac{\hw_j}{n\Opw_j}\Bigr)\Bigr\}$. This leads
immediately to the following corollary.  
\begin{corollary}Let the assumptions of Theorem~\ref{D:thm:upper:unknown} be
  satisfied. If in addition $(\hw/\Opw)$ is 
  non-decreasing and we have $\sup_{n\in\NN}\{R_n^\diamond / R_n^\ast\}<\infty$, then
  \begin{equation*}
 \sup_{P_U\in \cU_\sigma} \sup_{\sol \in\cF_\solw^{\solr}}
    \left\{
      \Ex\normV{\hsol_{\hk}-\sol}_\hw^2\right\}=O(\dstar),\qquad\mbox{ as }\quad n\to\infty.
  \end{equation*}
\end{corollary}
A sufficient condition for the two rates $R_n^\diamond$ and
$R_n^\ast$ to be of the same order is obviously given by
$\sup_{k\geq1}\{\delta_k/(\sum_{1\leq j\leq k}\hw_j/\Opw_j)\}<\infty$
and $\sup_{n\in\N}(\kstar/ \Cy_n^l)\leq 1$. However, this condition is
not necessary to establish the order optimality of the estimator as
follows from the example {\bf[is]} below.


\subsection{Illustration: estimation of derivatives  (continued)} 
The following result shows that even without any prior knowledge on
the structural function $\sol$ and for all smoothing operators $\Op$,
the fully adaptive penalized estimator automatically attains the
optimal rate in the finitely and in the infinitely smoothing case.
Recall that the computation of the dimension parameter $\hk$ given
in~\eqref{D:eq:21} involves the sequence $\Cy^u$, which in
our illustration satisfies $\Cy_n^u \asymp n^{1/(2s)}$ since
$\hw_{j} = j^{2s}$, $j\geq1$.
\begin{proposition}\label{D:prop:ex-cont-ada-ukn} Suppose that the
  i.i.d. $(Y,Z,W)$-sample of size $n$ obeys the model
  (\ref{D:model:NP}--\ref{D:model:NP2}) and that $P_U\in\cU_\sigma$,
  $\sigma>0$.  Consider the estimator $\hsol_{\hk}$ given in
  \eqref{D:gen:def:est} with~$\hk$ defined by~\eqref{D:eq:21}.
  \begin{itemize}
  \item[{\bf[fs]}] In the finitely smoothing case, we obtain\\[1ex]
    \hspace*{5ex}$ \sup_{P_U\in \cU_\sigma} \sup_{\sol \in\cW_p^\solr}     \left\{ \Ex\normV{\hsol^{(s)}_{\hk}-\sol^{(s)}}^2\right\}=
      O(n^{-2(p-s)/(2p+2a+1)})$.
  \item[\bf{[is]}] In the infinitely smoothing case, we have\\[1ex]
\hspace*{5ex}$ \sup_{P_U\in \cU_\sigma} \sup_{\sol \in\cW_p^\solr}
    \left\{ \Ex\normV{\hsol^{(s)}_{\hk}-\sol^{(s)}}^2\right\} =
    O((\log n)^{-(p-s)/a})$.
  \end{itemize}
\end{proposition}


\section{Conclusion}
\label{D:sec:conlusion}

In this work, we have developed a minimax theory for the estimation of
the structural function in a nonparametric regression model with
instrumental variables. We have defined a least squares estimator
involving dimension reduction and additional thresholding
which can attain the minimax optimal rate when the dimension parameter
is chosen in an appropriate way. This choice, however, depends on
characteristics of the structural function and the conditional expectation operator which are not
known in practice.

In order to solve this problem, we have proposed a data-driven
estimator which attains the minimax optimal rate over a wide range
of classes of structural functions. Unfortunately, we still need the additional assumption that
the eigenfunctions of the conditional expectation operator are known,
in which case the proposed estimator takes the form of an orthogonal
series. Furthermore, we show in~\eqref{D:eq:1} that if $\cS_k$
is the subspace generated by the first $k$ eigenfunctions, then we
have for all $k\leq k'$ and $t\in \cS_k$ that $\skalarV{ t,
  \hsol_{k'}}_\hw = \skalarV{ t, \hsol_{k}}_\hw$. If, however, $\cS_k$
is generated by an arbitrary set of linearly independent functions,
this is not true in general. In particular, the
estimate~\eqref{D:eq:4} on which the proof is essentially based, does
not hold anymore.

\cite{LM:09} develop, also under the assumption of known
eigenfunctions, oracle inequalities for an adaptive estimator which
attains the optimal bounds up to a logarithmic loss. In contrast to
this, the method we presented in this article does not suffer from
this loss.

\enlargethispage{2em}
The extension of this methodology to the case where the eigenfunctions
are unknown is an interesting problem worth investigating in the near
future.


\appendix

\section{Proofs}
\label{sec:proofs}

\subsection*{Minimax theory: lower bound}
\label{sec:minimax-theory}

\proofof{Theorem~\ref{D:res:lower}} Consider a pair $(Z,W)$ with
associated conditional expectation operator $T\in \cTdw$. Let
\[\zeta:=\kappa\min(\solr, 1/(2\Opd)) \text{\quad and\quad}  \alpha_n:=\dstar
(\sum_{j=1}^{\kstar} \hw_j/(\Opw_j n))^{-1}.\] Then,  the function
$\sol:= (\zeta\alpha_n/n)^{1/2}\sum_{j=1}^{\kstar}\Opw_j^{-1/2}
\basZ_j$ belongs to the class $\cF_\solw^\solr$, because the
monotonicity of $(\solw/\hw)$ implies
$\normV{\sol}^2_\solw\leq
\solr \kappa (\solw_{\kstar}/\hw_{\kstar}) \dstar\leq \solr$, using
successively the definitions of $\alpha_n$ and $\kappa$.  
Based on $\sol$, the candidates for the structural function are defined as
\[\sol_\theta:= \sum_{j=1}^{\kstar} \theta_j[\sol]_j \basZ_j\]
for every $\theta:=(\theta_j)\in\{-1,1\}^{\kstar}$. These functions
obviously belong to $\cF_\solw^\solr$, too.
Let $V\sim\cN(0,1)$ be a random variable independent of $(Z,W)$. For
every $\theta:=(\theta_j)\in\{-1,1\}^{\kstar}$, the distribution of
the random variable
\[U_\theta:= [T\sol_{\theta}](W)-\sol_{\theta}(Z) +V\]
then belongs to $\cU_\sigma$ for all $ \sigma^4\geq 8(3+2\solr^2
\Gamma^2\eta)$: Firstly, $\Ex[U_\theta|W=0]$. Secondly, 
we have
\[|\Ex
[f(Z)|W]|^4\leq 
\solr^2 \Gamma \sum_{j\in\N}\solw_j^{-1}\Ex[\basZ^4_j(Z)|W]\leq
\solr^2\Gamma^2\eta\]
for all   all $f\in \cF_\solw^\solr$, which follows from  
 the condition
$\Gamma=\sum_{j\in\N}\solw_j^{-1}<\infty$ together with
$\sup_j\Ex[\basZ_j^4(Z)|W]\leq \eta$, applying the
Cauchy-Schwarz inequality twice.  From this
estimate we conclude $\Ex [\sol_\theta^4(Z)|W]\leq
\eta\solr^2\Gamma^2$ and $|[T\sol_\theta](W)|^4\leq \Ex
[\sol_\theta^4(Z)|W]\leq \eta\solr^2\Gamma^2$. By combination of the
last two bounds we obtain $\Ex [U_\theta^4|W]\leq 8\{
2\eta\solr^2\Gamma^2 +3\}$.

Consequently, for any $\theta$, the tuple $(Y,Z,W)$ defined by
$Y:=\sol_{\theta}(Z)+U_\theta$ obeys the model
(\ref{D:model:NP}--\ref{D:model:NP2}). Let
$(Y_i,Z_i,W_i)_{i=1,\ldots,n}$ be $n$ iid.\ copies of $(Y,Z,W)$ and
denote their joint distribution by $P_{\theta}$.

Under the law $P_\theta$, the conditional
distribution of $Y_i$ given $W_i$ is then Gaussian with mean $
[T\sol_\theta](W_i)$ and variance $1$.  Furthermore, for
$j=1,\dotsc,\kstar$ and for each $\theta$ we introduce $\theta^{(j)}$ by
$\theta^{(j)}_{l}=\theta_{l}$ for $j\ne l$ and
$\theta^{(j)}_{j}=-\theta_{j}$. Then, it is easily seen that the
log-likelihood of ${P}_{\theta}$ with respect to ${P}_{\theta^{(j)}}$
is given by
\begin{equation*}
  \log\Bigl(\frac{d{P}_{\theta}}{d{P}_{\theta^{(j)}}}\Bigr)=\sum_{i=1}^n 2(Y_i - [T\sol_\theta](W_i)) \theta_j [\sol]_j[T\basZ_j](W_i) + 2[\sol]_j^2\sum_{i=1}^n  |[T\basZ_j](W_i)|^2.
\end{equation*}
Its expectation with respect to ${P}_{\theta}$ satisfies
\[\Ex_{{P}_{\theta}}[\log(d{P}_{\theta}/d{P}_{\theta^{(j)}})]=
2n[\sol]_j^2 \normV{T\basZ_j}_W^2 \leq 2n \Opd [\sol]_{j}^2 \Opw_{j},\]
because $T\in\cTdw$. \label{sym:kull1}
In terms of the Kullback-Leibler divergence,
this means \[KL(P_{\theta},P_{\theta^{(j)}})\leqslant 2\,\Opd\,n\,
[\sol]_{j}^2 \Opw_{j}.\] Since the Hellinger distance satisfies
$H^2(P_{\theta},P_{\theta^{(j)}}) \leqslant
KL(P_{\theta},P_{\theta^{(j)}})$, we can use
the definition of $\sol$, the property $\alpha_n\leq \kappa^{-1}$, and
the definition of $\zeta$ successively and obtain that
\begin{equation}\label{D:pr:lower:e3}
  H^2(P_{\theta},P_{\theta^{(j)}})  \leqslant 2\,\Opd\,n \, [\sol]_{j}^2 \Opw_{j}\leq 2\Opd \, \zeta\, \alpha_n \leq  1.
\end{equation} 
Considering the Hellinger affinity $\rho(P_{\theta},P_{\theta^{(j)}})=
\int \sqrt{dP_{\theta}dP_{\theta^{(j)}}}$, we can write for any
estimator $\tsol$ of $\sol$ that
\begin{align*}
  \rho(P_{\theta},P_{\theta^{(j)}})&\leqslant \int
  \frac{|[\tsol-\sol_{\theta^{(j)}}]_j|}{|[\sol_\theta-\sol_{\theta^{(j)}}]_j|}\sqrt{dP_{\theta}dP_{\theta^{(j)}}}
  +
  \int\frac{|[\tsol-\sol_{\theta}]_j|}{|[\sol_\theta-\sol_{\theta^{(j)}}]_j|}
  \sqrt{dP_{\theta}dP_{\theta^{(j)}}}\\
  &\leqslant \Bigl( \int
  \frac{|[\tsol-\sol_{\theta^{(j)}}]_j|^2}{|[\sol_\theta-\sol_{\theta^{(j)}}]_j|^2}
  dP_{\theta^{(j)}} \Bigr)^{1/2} + \Bigl( \int
  \frac{|[\tsol-\sol_{\theta}]_j|^2}{|[\sol_\theta-\sol_{\theta^{(j)}}]_j|^2}
  dP_{\theta}\Bigr)^{1/2}.
\end{align*}
Rewriting the last estimate using the identity
$\rho(P_{\theta},P_{\theta^{(j)}})=1-\frac{1}{2}H^2(P_{\theta},P_{\theta^{(j)}})$
and~\eqref{D:pr:lower:e3}, we obtain
\begin{equation*}
  \Bigl\{\Ex_{{\theta}}|[\tsol-\sol_\theta]_j|^2+\Ex_{{\theta^{(j)}}}|[\tsol-\sol_{\theta^{(j)}}]_j|^2 \Bigr\} \geq  \frac{1}{8} |[\sol_\theta-\sol_{\theta^{(j)}}]_j|^2=\frac{1}{2}[\sol]_j^2. \end{equation*}
We combine the last estimate with  the following reduction scheme,
which  is the key argument of this proof: 
\begin{align*}
  \sup_{P_U\in\cU_\sigma} \sup_{\sol \in \cF_\solw^\solr}
  &\Ex_{P_\theta}\normV{\tsol -\sol}_\hw^2 \geq \sup_{\theta\in
    \{-1,1\}^{\kstar}} \Ex_{P_\theta}\normV{\tsol
    -\sol_\theta}_\hw^2\\
  &\geq \frac{1}{2^{{\kstar}}}\sum_{\theta\in \{-1,1\}^{\kstar}}\sum_{j=1}^{\kstar}\hw_j\Ex_{P_{\theta}}|[\tsol-\sol_\theta]_j|^2\\
  &= \frac{1}{2^{{\kstar}}}\sum_{\theta\in
    \{-1,1\}^{\kstar}}\sum_{j=1}^{\kstar}\frac{\hw_j}{2}\Bigl\{\Ex_{P_{\theta}}|[\tsol-\sol_\theta]_j|^2+\Ex_{P_{\theta^{(j)}}}|[\tsol-\sol_{\theta^{(j)}}]_j|^2
  \Bigr\}\\
  &\geq \frac{1}{2^{{\kstar}}}\sum_{\theta\in
    \{-1,1\}^{\kstar}}\sum_{j=1}^{\kstar}\frac{\hw_j}{4}[\sol]_j^2=\frac{\zeta\alpha_n}{4}\sum_{j=1}^{\kstar}\frac{\hw_j}{n\Opw_j}.
\end{align*}
Hence, from the definition of $\zeta$ and $\alpha_n$ we obtain the
lower bound given in the theorem.\qed


\subsection*{Minimax theory: upper bounds}\label{D:proofs:upper:bounds:mmax}
We begin by defining and recalling notations to be used without
further reference in the proofs of this section. Given $k>0$, denote
$\sol_k:=\sum_{j=1}^k\fou{\sol_k}_{j}\basZ_j$ with
$\fou{\sol_k}_{\uk}=\fou{\Op}_{\uk}^{-1}\fou{\gf}_{\uk}$ which is
well-defined since $\fou{\Op}_{\uk}$ is non singular. Then, the
identities $[\Op(\sol-\sol_k)]_{\uk}=0$ and $[\sol_k-E_k \sol]_{\uk} =
\fou{\Op}_{\uk}^{-1}\fou{\Op E_k^\perp \sol}_{\uk}$ hold true.
Furthermore, let $\fou{\Xi}_{\uk}:= \hfou{\Op}_{\uk}- \fou{\Op}_{\uk}$
and define vectors $\fou{B}_{\uk}$ and $\fou{S}_{\uk}$ by
\begin{gather*}
  \fou{B}_j:=\frac{1}{n}\sum_{i=1}^n U_i \basW_j(W_i)\qquad\text{and}\qquad
  \fou{S}_j:=\frac{1}{n}\sum_{i=1}^n \basW_j(W_i)\{ \sol(Z_i) -   
\fou{\sol_k}_{\uk}^t\fou{\basZ(Z_i)}_{\uk}\},\;1\leq j\leq k,
\end{gather*} 
such that $\hfou{\gf}_{\uk}- \hfou{\Op}_{\uk}
\fou{\sol_k}_{\uk}=\fou{B}_{\uk}+ \fou{S}_{\uk}$. Note that $ \Ex
\fou{B}_{\uk}=0$ due to the mean independence, i.e., $\Ex[U|W]=0$, and
that $\Ex\fou{S}_{\uk}= \fou{\Op\sol}_{\uk} - \fou{\Op\sol_k}_{\uk}=
0$. Moreover, let us introduce the events
\begin{equation*}\label{def:OmOm}
  \Omega:=\{ \normV{\hfou{\Op}^{-1}_{\uk}}\leq \sqrt n\}\quad\text{
    and }\quad
  \Omega_{1/2}:=
  \{\normV{\fou{\Xi}_{\uk}}\,\normV{\fou{\Op}_{\uk}^{-1}}\leq 1/2\}.
\end{equation*}
At the end of this section we shall prove some technical lemmas
(Lemmas~\ref{D:app:gen:upper:l2}-- \ref{D:app:gen:upper:l4}) which are used
in the following proofs.
\bigskip

\proofof{Theorem \ref{D:res:upper:A3}} Define
$\tsol_{\kstar}:=\sol_{\kstar}\1_\Omega$ and decompose the risk into
two terms,
\begin{equation}\label{D:app:gen:dec}
  \Ex\normV{\hsol_{\kstar}-\sol}^2_\hw \leq 2 \{ \Ex\normV{\hsol_{\kstar}-\tsol_{\kstar}}^2_\hw+ \Ex\normV{\tsol_{\kstar}-\sol}^2_\hw\}=:2\{A_1+A_2\},
\end{equation}
which we bound separately. Consider first $A_2$. Using the identity
$\normV{\tsol_{\kstar}-\sol}^2_\hw=
\normV{\sol_{\kstar}-\sol}^2_\hw\1_\Omega
+\normV{\sol}^2_\hw\1_{\Omega^c}$ we deduce
\[\Ex\normV{\tsol_{\kstar}-\sol}^2_\hw \leq
\normV{\sol_{\kstar}-\sol}^2_\hw + \normV{\sol}^2_\hw
P(\Omega^c).\] 
Since $(\hw/\solw)$ is monotonically decreasing,
the last estimate together with \eqref{D:app:gen:upper:l3:e1} in
Lemma~\ref{D:app:gen:upper:l3} implies for all
$\sol\in\cF_\solw^\solr$
\begin{align}  \label{D:eq:19}
  \begin{split}
    \Ex\normV{\tsol_{\kstar}-\sol}^2_\hw &\leq 4\,\OpD\, \Opd \,
    \solr\, \dstar\,\max\bigg( 1, \frac{\Opw_{\kstar}}{\hw_{\kstar}}
    \max_{1\leq j\leq \kstar}\frac{\hw_j}{\Opw_j}\bigg) + P(\Omega^c)\\
    &\leq \dstar\;\bigg\{4\,\OpD\, \Opd \, \solr\,\max\left( 1,
      \frac{\Opw_{\kstar}}{\hw_{\kstar}} \max_{1\leq j\leq
        \kstar}\frac{\hw_j}{\Opw_j}\right) + C(\solw,\Opw,\eta,\OpD)
    \bigg\}
  \end{split}
\end{align}
by employing the definition of $\dstar$ and applying
Lemma~\ref{lem:a4fromcb}.

\noindent Consider $A_1$.  From the identity
$\hfou{\gf}_{\ukstar}-\hfou{\Op}_{\ukstar}\fou{\sol_m}_{\ukstar}=\fou{B}_{\ukstar}
+ \fou{S}_{\ukstar}$ follows
\begin{align*}
  \fou{\hsol_{\kstar}-&\tsol_{\kstar}}_{\ukstar}=\{\fou{\Op}_{\ukstar}^{-1}+\fou{\Op}_{\ukstar}^{-1}(\fou{\Op}_{\ukstar} - \hfou{\Op}_{\ukstar})\hfou{\Op}_{\ukstar}^{-1}\}\{\fou{B}_{\ukstar} + \fou{S}_{\ukstar}\}\1_\Omega \\
  &= \fou{\Op}_{\ukstar}^{-1} \,\{\fou{B}_{\ukstar} +
  \fou{S}_{\ukstar}\}\1_{\Omega}-\fou{\Op}_{\ukstar}^{-1}\fou{\Xi}_{\ukstar}\hfou{\Op}_{\ukstar}^{-1}\,\{\fou{B}_{\ukstar}
  + \fou{S}_{\ukstar}\}\1_\Omega.
\end{align*}
By making use of this identity we decompose $A_1$ further into two
terms
\begin{multline}\label{D:eq:16}
  \Ex\normV{\hsol_{\kstar}-\tsol_{\kstar}}^2_\hw \leq 2 \Ex
  [\normV{\fou{\diag(\hw)}_{\kstar}^{1/2}\fou{\Op}_{\kstar}^{-1}
    \,\{\fou{B}_{\kstar} + \fou{S}_{\kstar}\} }^2\1_{\Omega}]
  \\\hfill+ 2 \Ex
  [\normV{\fou{\diag(\hw)}_{\kstar}^{1/2}\fou{\Op}_{\kstar}^{-1}\fou{\Xi}_{\kstar}\hfou{\Op}_{\kstar}^{-1}\,\{\fou{B}_{\kstar}
    + \fou{S}_{\kstar}\}}^2\1_\Omega]=:2\{A_{11}+A_{12}\}
\end{multline}
which we bound separately. In case of $A_{11}$ we employ
successively~\eqref{D:app:gen:upper:l2:e4} from
Lemma~\ref{D:app:gen:upper:l2} with
$M:=\fou{\diag(\hw)}_{\kstar}^{1/2}\fou{\Op}_{\kstar}^{-1}$, the
elementary inequality $\tr(A^tB^tBA)\leq \normV{A}^2 \tr(B^tB)$ valid
for all $(k\times k)$ matrices $A$ and $B$ and the extended link
condition \eqref{D:bm:link:gen}, that is,
$\normV{\fou{\diag(\Opw)}_{\kstar}^{1/2}\fou{\Op}_{\kstar}^{-1}}^2\leq
\OpD$. Thereby, we obtain
\begin{multline}\label{D:app:gen:upper:e1:1}
  \Ex
  [\normV{\fou{\diag(\hw)}_{\ukstar}^{1/2}\fou{\Op}_{\ukstar}^{-1}
    \,\{\fou{B}_{\ukstar} + \fou{S}_{\ukstar}\} }^2\1_{\Omega}] \\\leq
  (2/n)\, \OpD\,
  \tr\bigg(\fou{\diag(\Opw)}_{\ukstar}^{-1/2}\fou{\diag(\hw)}_{\ukstar}\fou{\diag(\Opw)}_{\ukstar}^{-1/2}\bigg)
  \{\sigma^2+\eta^2\,\Gamma\, \normV{ \sol -
    \sol_{\kstar}}_\bw^2\}\\
  = 2\OpD \{\sigma^2+\eta^2\,\Gamma\, \normV{ \sol -
    \sol_{\kstar}}_\bw^2\}\sum_{j=1}^{\kstar}\frac{\hw_j}{n\Opw_j}.
\end{multline}
Consider now $A_{12}$. Observe that $\normV{\fou{\diag(\hw)}_{\ukstar}^{1/2}\fou{\Op}_{\ukstar}^{-1}}^2\leq
\OpD \max_{1\leq j\leq\kstar}{\hw_j}/{\Opw_j}$ for all $T\in\cTdDw$.
Applying the last inequality together with
\[\normV{\hfou{\Op}_{\ukstar}^{-1}}^2\1_{\Omega_{1/2}} \leq 4
\OpD/\Opw_{\kstar} \quad\text{ and }\quad \normV{\hfou{\Op}_{\ukstar}^{-1}}^2\1_\Omega
\leq n,\] we see that there exists a numerical constant $C>0$ such that
\begin{align*}
  \Ex
  [\normV{\fou{\diag(\hw)}_{\ukstar}^{1/2}&\fou{\Op}_{\ukstar}^{-1}\fou{\Xi}_{\ukstar}\hfou{\Op}_{\ukstar}^{-1}\,\{\fou{B}_{\ukstar}
    + \fou{S}_{\ukstar}\}}^2\1_\Omega]\\&\leqslant \OpD \max_{1\leq
    j\leq\kstar}\frac{\hw_j}{\Opw_j} \Bigl\{ 4 \OpD\Opw_{\kstar}^{-1}
  \Ex\normV{ \fou{\Xi}_{\ukstar}}^2\normV{\fou{B}_{\ukstar}
    +\fou{S}_{\ukstar}}^2\1_{\Omega_{1/2}} 
\\ & \hspace{10em}+
  n \Ex\normV{ \fou{\Xi}_{\ukstar}}^2\normV{\fou{B}_{\ukstar} +\fou{S}_{\ukstar}}^2\1_{\Omega_{1/2}^c}\Bigr\}\\
  &\leq \OpD \max_{1\leq j\leq\kstar}\frac{\hw_j}{\Opw_j} \Bigl\{ 4
  \OpD\Opw_{\kstar}^{-1} \bigl(\Ex
  \normV{\fou{\Xi}_{\ukstar}}^4\bigr)^{1/2} \\& \hspace{3em} + n \bigl(\Ex
  \normV{\fou{\Xi}_{\ukstar}}^8\bigr)^{1/4} P(\Omega_{1/2}^c)^{1/4}\Bigr\}\bigl(\Ex\normV{\fou{B}_{\ukstar}+\fou{S}_{\ukstar}}^4\bigr)^{1/2}\\
 & \leq C\,\max_{1\leq j\leq\kstar}\frac{\hw_j}{n\Opw_j}
  \,\OpD\,\eta^4\, ( \sigma^2+\Gamma\, \normV{ \sol - \sol_{\kstar}}_\bw^2 )
\\ &  \hspace{10em} \Bigl\{ 4\OpD \frac{ (\kstar)^3}{\Opw_{\kstar} n} + (\kstar)^3
  |P(\Omega^c_{1/2})|^{1/4} \Bigr\}
\end{align*}
where the last bound follows from \eqref{D:app:gen:upper:l2:e2:1},
\eqref{D:app:gen:upper:l2:e2:2} and \eqref{D:app:gen:upper:l2:e3} in Lemma~\ref{D:app:gen:upper:l2}.  By combination of the last bound and
\eqref{D:app:gen:upper:e1:1} via the decomposition \eqref{D:eq:16} there
exists a numerical constant $C>0$ such that
\begin{gather*}
  \Ex\normV{\hsol_{\kstar}-\tsol_{\kstar}}^2_\hw \leq C  \,\OpD\,\eta^4\, ( \sigma^2+\Gamma\, \normV{ \sol -   \sol_{\kstar}}_\bw^2 ) 
  \Bigl\{ 4\OpD \zeta/\kappa + (\kstar)^3 |P(\Omega^c_{1/2})|^{1/4} \Bigr\} \sum_{j=1}^{\kstar}\frac{\hw_j}{n\Opw_j} .
\end{gather*}
Furthermore, taking into account the estimate
\eqref{D:app:gen:upper:l3:e1} from Lemma~\ref{D:app:gen:upper:l3} with
$\hw=\solw$ and the definition of $\dstar$, the last inequality implies
\begin{equation*}
  \Ex\normV{\hsol_{\kstar}-\tsol_{\kstar}}^2_\hw \leq C  \,\OpD\,\eta^4\, ( \sigma^2+4\Gamma\, \OpD\Opd\solr ) 
  \Bigl\{ 4\OpD \zeta/\kappa + (\kstar)^3 |P(\Omega^c_{1/2})|^{1/4} \Bigr\} \dstar.
\end{equation*}
Finally, using the decomposition~\eqref{D:app:gen:dec}, the result of
the theorem follows from the last estimate and \eqref{D:eq:19}, since
$(\kstar)^3(P(\Omega_{1/2}))^{1/4}\leq C(\solw,\Opw,\eta,\OpD)$ by Lemma~\ref{lem:a4fromcb}.\qed

\subsection*{Illustration}
\label{sec:illustration}

\proofof{Proposition~\ref{D:coro:ex:lower}} Since for each $0\leqslant
s \leqslant p$ we have $\Ex\normV{\txdf^{(s)}-\xdf^{(s)}}^2 \sim \Ex
\normV{\txdf-\xdf}_{\hw}^2$ we apply the general result given
Theorem~\ref{D:res:lower}. In both cases, the additional conditions
formulated in Theorem~\ref{D:res:lower} are easily
verified. Therefore, it is sufficient to evaluate the lower bound
$\dstar$ given in \eqref{D:def:kn}. Note that the optimal dimension
parameter $\kstar$ satisfies $\dstar\sim\hw_{\kstar}/\solw_{\kstar}
\sim \sum_{l=1}^{\kstar}\hw_l/(n\Opw_l)$ since both sequences
$(\xdfw_j/\hw_j)$ and $(\sum_{0<|l|\leq j}\frac{\hw_l}{n\edfw_l})$ are
non-increasing.

{\bf[p-$\pmb{\Opw}$]} The well-known approximation $\sum_{j=1}^{k}
j^{r}\sim k^{r+1}$ for $r>0$ implies\\
$n\sim(\solw_{\kstar}/\hw_{\kstar})\sum_{l=1}^{\kstar}\hw_l/\Opw_l
\sim (\kstar)^{2a+2p+1}$. It follows that $\kstar\sim n^{1/(2p+2a+1)}$
and the lower bound writes $\dstar \sim n^{-(2p-2s)/(2p+2a+1)}$.

{\bf[e-$\pmb{\Opw}$]} Applying Laplace's Method (c.f.\ Chapter 3.7 in
\cite{Olver1974}) we have\\
$n\sim(\solw_{\kstar}/\hw_{\kstar})\sum_{l=1}^{\kstar}\hw_l/\Opw_l\sim
(\kstar)^{2p}\exp(|\kstar|^{2a})$ which implies that\\ $\kstar\sim
\{\log (n/ (\log n)^{p/a})\}^{1/(2a)} = (\log n)^{1/(2a)}(1+o(1))$ and
that the lower bound can be rewritten as $\dstar\sim (\log
n)^{-(p-s)/a}$.  \qed

\proofof{Proposition~\ref{D:coro:ex:upper}} 
Since in both cases the dimension parameter is chosen in the optimal way (see the
proof of Proposition~\ref{D:coro:ex:lower}), the result follows from
Theorem~\ref{D:res:upper:A3}.\qed


\subsection*{Adaptive estimation}
\label{sec:adaptive-estimation}
Before proving Theorem~\ref{D:thm:upper:unknown}, we define some
notation to be used in the proof. Given $u\in L^2[0,1]$ we denote by
$[u]$ the infinite vector of Fourier coefficients
$[u]_j:=\skalarV{u,\basTr_j}$. In particular we use the notations
\begin{multline*}
\hspace{-1em}  \hsol_k=
  \sum_{j=1}^k\frac{\hfou{\gf}_j}{\hfou{\Op}_{jj}}\1\{\inf_{1\leq
    j\leq k}\hfou{\Op}_{jj}^2\geq 1/n\} \basTr_j, 
\; \tsol_k:=  \sum_{j=1}^k\frac{\hfou{\gf}_j }{\fou{\Op}_{jj}}\basZ_j,
 \; \sol_k:=\sum_{j=1}^k\frac{\fou{\gf}_j}{
    \fou{\Op}_{jj}} \basTr_j, \\
  \hPhi_{u}:=
  \sum_{j\in\N}\frac{\fou{u}_j}{\hfou{\Op}_{jj}}\1\{\inf_{1\leq j\leq
    k}\hfou{\Op}_{jj}^2\geq 1/n\}\basTr_j,\quad \tPhi_{u}:=
  \sum_{j\in\N}\frac{\fou{u}_j}{\fou{\Op}_{jj}}\basTr_j.
  \hfill\end{multline*}

\noindent Furthermore, let $\hgf$ be the function with Fourier coefficients
$[\hgf]_j:=\hfou{\gf}_j$ and observe that $\Ex \hgf=\gf$.  Given
$1\leq k\leq k'$ we have then for all $t\in\cS_k:={\rm
  span}\{\basTr_{1},\dotsc,\basTr_k\}$
 \begin{align}\label{D:eq:1}\notag
\skalarV{ t, \hsol_{k'}}_\hw & =  \skalarV{ t,
\hPhi_{\hgf}}_\hw=\frac{1}{n} \sum_{i=1}^n  \sum_{j=1}^k Y_i\basTr_j(W_i) 
\frac{\hw_j[t]_j}{\hfou{\Op}_{jj}}\I{[\inf_{1\leq j\leq k}\hfou{\Op}_{jj}^2\geq 1/n]}=\skalarV{ t, \hsol_{k}}_\hw,\\
 \skalarV{ t,  \tsol_{k'}}_\hw&=\skalarV{ t,\tPhi_{\hydf}}_\hw=  \frac{1}{n} \sum_{i=1}^n  \sum_{j=1}^k Y_i\basTr_j(W_i)
\frac{\hw_j[t]_j}{\fou{\Op}_{jj}}= \skalarV{ t,  \tsol_{k}}_\hw,
\hfill\\\notag
\skalarV{ t, \sol_{k'}}_\hw&= \skalarV{ t,
\tPhi_{\gf}}_\hw= \sum_{j=1}^k\frac{\hw_j[t]_j[\gf]_j}{\fou{\Op}_{jj}}=  \sum_{j=1}^k \hw_j[t]_j[\sol]_j=\skalarV{ t, \sol}_\hw.\hfill
\end{align}
  Consider the contrast
$\ct(t) :=   \normV{t}^2_\hw - 2 {\skalarV{ t, \hPhi_{\hgf}}}_\hw$, for all $t\in L^2[0,1]$. Obviously it follows for all $t\in\cS_k$  that 
$\ct(t) = \normV{t-\hsol_k}^2_\hw- \normV{\hsol_k}^2_\hw$ and, hence 
\begin{equation}\label{D:eq:2}
\arg\min_{t\in\cS_k}  \ct(t) = \hsol_{k},\quad \forall \,k\geq 1.
\end{equation}
In order to decompose the risk, we will use the events
\begin{align}
  \label{eq:OmOmpq}
  \begin{split}
    \Omega_q &:= \bigg\{ \forall\; 1\leq j\leq \Cy_n\; \bigg|\;
    \Big|\hfou{\Op}_{jj}^{-1}-\fou{\Op}_{jj}^{-1}\Big| \leq
    \frac{1}{2|\fou{\Op}_{jj}|} \;\wedge\;
    \hfou{\Op}_{jj}^2 \geq 1/n     \bigg\} \\
    \Omega_{p}&:=\Bigl\{ \Cy^l_n \leq \hCy_n \leq \Cy_n\Bigr\}
  \end{split}
\end{align}
in the following proofs. 

\bigskip
\proofof{Theorem~\ref{D:thm:upper:unknown}} 
Given the preliminary considerations above, the adaptive choice $\whk$ of the dimension parameter can be rewritten as
\begin{equation}\label{D:eq:3}
\whk = \argmin_{1\leq k\leq \hCy_n}\left\{\ct(\hsol_{k}) + \hpen(k)  \right\}\hspace*{5ex}\mbox{with}\hspace*{5ex}\hpen(k):=  540\,\Ex[Y^2]\frac{ \widehat\delta_k}{n} .
\end{equation}
Let us suppose in the following paragraph that the event $\Omega_p$ occurs. Then we have for all
$1\leq k \leq \Cy^l_n$ that
\[\ct(\hsol_{\hk}) + \hpen(\hk) \leq \ct(\hsol_{k}) + \hpen(k) \leq  \ct(\sol_k) + \hpen(k),\]
using first~\eqref{D:eq:3} and then~\eqref{D:eq:2}.  This inequality
implies
\begin{align*}
  \normV{\hsol_{\hk}}^2_\hw-\normV{\sol_k}^2_\hw
  &\leq 2\skalarV{\hsol_{\hk}-\sol_k, \hsol_{\hk} }_\hw
  + \hpen(k)-\hpen(\hk),  \end{align*} 
which together with the identities given in \eqref{D:eq:1} implies for all $1\leq k \leq \Cy^l_n$ 
\begin{multline}\label{D:eq:4}
    \normV{\hsol_{\hk}-\sol}^2_\hw =  \normV{\sol-\sol_{k}}^2_\hw +  \normV{\hsol_{\hk}}^2_\hw-\normV{\sol_k}^2_\hw -  2\skalarV{\hsol_{\hk} - \sol_k,\sol}_\hw\\
    \hfill\leq \normV{\sol-\sol_{k}}^2_\hw + \hpen(k)-\hpen(\hk) 
+2\skalarV{\hsol_{\hk} - \sol_k,\hPhi_{\hgf}-\tPhi_{\gf}}_\hw
\end{multline}
Consider the unit ball $\cB_k:=\{f\in \cS_k\;|\;\normV{f}_\hw\leq 1\}$ and,   for arbitrary $\tau>0$ and $t\in \cS_k$, the elementary
inequality 
\begin{gather*}
  2|\skalarV{t,h}_\hw|
\leq
  2\normV{t}_\hw \sup_{t\in \cB_k}|\skalarV{t,h}_\hw|
 \leq \tau
  \normV{t}_\hw^2+\frac{1}{\tau} \sup_{t\in \cB_k}|\skalarV{t,h}_\hw|^2= \tau
  \normV{t}_\hw^2+ \frac{1}{\tau}\sum_{j=1}^k\hw_j|\fou{h}_j|^2.
\end{gather*}
Combining the last estimate with (\ref{D:eq:4}) and $\hsol_{\whk}-\sol_k \in \cS_{\whk\vee k}\subset\cS_{\Cy_n} $   we obtain
\begin{multline*}
  \normV{\hsol_{\hk} -\sol}^2_\hw \leq \normV{\sol-\sol_k}^2_\hw
  + \tau\, \normV{\hsol_{\hk} -\sol_k}^2_\hw + \hpen(k) - \hpen(\hk)
   + \frac{1}{\tau}\supt{\Cy_n}|\skalarV{t,\hPhi_{\hgf}-\tPhi_\gf}_\hw|^2.
\end{multline*}
Letting $\tau := {1}/{3}$  it follows from $\normV{\hsol_{\hk}-\sol_k}^2_\hw \leq 2\normV{\hsol_{\hk}-\sol}^2_\hw
+2\normV{\sol_k-\sol}^2_\hw$ that
\begin{equation*}
 \frac{1}{3} \normV{\hsol_{\hk} -\sol}^2_\hw \leq \frac{5}{3}\normV{\sol-\sol_k}^2_\hw + \hpen(k) - \hpen(\hk)
   + 3\supt{\Cy_n}|\skalarV{t,\hPhi_{\hgf}-\tPhi_\gf}_\hw|^2.
\end{equation*}
Consider the functions $\hnu$ and $\hmu$  with  
\[\fou{\hnu}_j = \frac{1}{n}\sum_{i=1}^n
Y_i\I{[|Y_i|\leq\an]}\basTr_j(W_i) \quad\text{and}\quad
\fou{\hmu}_j = \frac{1}{n}\sum_{i=1}^nY_i\I{[|Y_i|>\an]}\basTr_j(W_i)\] respectively,
 as well as their centered versions  $\nu=\hnu-\Ex[\hnu]$ and $\mu=\hmu-\Ex[\hmu]$,  then we have
$\hgf-\gf=\nu+\mu$ and  
\begin{multline*}
 \frac{1}{3} \normV{\hsol_{\hk} -\sol}^2_\hw \leq \frac{5}{3}\normV{\sol-\sol_k}^2_\hw + \hpen(k) - \hpen(\hk)\\
\hfill   + 6\supt{\Cy_n}|\skalarV{t,\tPhi_{\nu}}_\hw|^2+12\supt{\Cy_n}|\skalarV{t,\hPhi_{\nu} - \tPhi_{\nu}}_\hw|^2+12\supt{\Cy_n}|\skalarV{t,\hPhi_{\mu}+\hPhi_{\ydf} - \tPhi_{\ydf}}_\hw|^2
\end{multline*}
Consider the decomposition  $|\skalarV{t,\hPhi_{\nu} -
    \tPhi_{\nu}}_\hw|^2 = |\skalarV{t,\hPhi_{\nu} -
    \tPhi_{\nu}}_\hw|^2 \I{\Omega_q} +|\skalarV{t,\hPhi_{\nu} -
    \tPhi_{\nu}}_\hw|^2 \I{\Omega_q^c}$.
Since $\I{[\hfou{\Op}_{jj}^2\geq1/n]}\I{\Omega_q} = \I{\Omega_q}$,   it follows that for all $1\leq j \leq \Cy_n$ we have
\begin{align*}
  \bigg(\frac{\fou{\Op}_{jj}}{\hfou{\Op}_{jj}}\I{[\hfou{\Op}_{jj}^2\geq1/n]} - 1\bigg)^2\,\I{\Omega_q} =  |\fou{\Op}_{jj}|^2\;\I{\Omega_q}\,
\bigg|\hfou{\Op}_{jj}^{-1}-\fou{\Op}_{jj}^{-1}\bigg|^2
\leq \frac{1}{4}.
\end{align*}
Hence,  
$  \supt{k} |\skalarV{t,\hPhi_{\nu} -
    \tPhi_{\nu}}_\hw|^2\,\I{\Omega_q} 
\leq \frac{1}{4} \supt{k} |\skalarV{t,\tPhi_\nu}_\hw|^2$ 
   for all $1\leq k\leq\Cy_n$ and
\begin{multline}\label{D:eq:6}
    \frac{1}{3} \normV{\hsol_{\hk} -\sol}^2_\hw \leq
    \frac{5}{3}\normV{\sol-\sol_k}^2_\hw +  9 \supt{k\vee \hk}
    |\skalarV{t,\tPhi_\nu}_\hw|^2  
    + \hpen(k) - \hpen(\hk)\hfill\\[.75em]\hfill
 +      12\supt{k\vee\hk}|\skalarV{t,\hPhi_{\nu} -
    \tPhi_{\nu}}_\hw|^2\,\I{\Omega_q^c}
+12\supt{k\vee\hk}|\skalarV{t,\hPhi_{\mu}+\hPhi_{\ydf} - \tPhi_{\ydf}}_\hw|^2.
\end{multline}
Define $\Delta^\Op_k:=\max_{1\leq
  j\leq k}\hw_j/|\fou{\Op}_{jj}|^2$, $\tau_k^\Op:=\max_{1\leq j\leq
  k}(\hw_j)_{\vee1}/|\fou{\Op}_{jj}|^2$,\\ and
  $\delta_k^\Op := {k}
    \Delta_k^\Op\;\big\{{\log (\tau_k^\Op\vee
      (k+2))}\;/\;{\log(k+2)}\big\}$. Then, it is easily seen that   
    \begin{equation} 
  \delta_k^\Op \leq \delta_k \, \Opd\, \frac{\log (3\Opd)}{\log
      3}=\delta_k\,\Opd\,\zeta_d  \qquad  \forall\;k\geq 1.\label{D:eq:7}
  \end{equation}
  with $\zeta_d= {(\log 3\Opd)}/{(\log 3)}$. Moreover, define the
  event $\Omega_{qp}:=\Omega_q\cap\Omega_{p}$ with~$\Omega_q$
  and~$\Omega_p$ from~\eqref{eq:OmOmpq}.
Observe that on $\Omega_q$ we have  $ (1/2)\Delta^\Op_k\leq \widehat\Delta_k\leq(3/2)\Delta^\Op_k$ for all $1\leq k\leq
\Cy_n$
and hence $(1/2)[\Delta^\Op_k\vee(k+2)]\leq[ \widehat\Delta_k\vee
(k+2)]\leq(3/2)[\Delta^\Op_k\vee(k+2)]$, which implies
\begin{align*}
  \begin{split}
    (1&/2) k \Delta^\Op_k\Bigl(\frac{ \log[ \Delta^\Op_k\vee
      (k+2)]}{\log(k+2)}\Bigr)\Bigl(1-\frac{\log 2}{\log
      (k+2)}\frac{\log (k+2)}{\log (\Delta^\Op_k \vee [k+2])}\Bigr) 
\\[1ex]&\leq
    \widehat\delta_k \leq (3/2) k \Delta^\Op_k\Bigl(\frac{\log (\Delta^\Op_k
      \vee [k+2])}{\log( k+2)}\Bigr) \Bigl(1+\frac{\log3/2}{\log(k+2)
    } \frac{\log (k+2)}{\log (\Delta^\Op_k \vee [k+2])}\Bigr).
  \end{split}
\end{align*}
Using $ {\log (\Delta^\Op_k
  \vee (k+2))}/{\log( k+2)}\geq 1$, we conclude from the last estimate that
\begin{align*}
  \begin{split}
{ \delta^\Op_k/10\leq}
(\log 3/2)/(2 \log 3) \delta^\Op_k&\leq (1/2) \delta^\Op_k[1-(\log
2)/\log(k+2)] \leq {\widehat\delta_k} \\ &\leq (3/2) \delta^\Op_k[1+
(\log3/2)/\log (k+2)]\leq { 3\delta^\Op_k}.
\end{split}
\end{align*} 
Recalling that $\hpen(k)= 540\,\Ex[Y^2]\,
\widehat\delta_k n^{-1}$, we define 
\begin{equation*}
  \label{D:eq:14}
\pen(k):= 54\,\Ex[Y^2]\,\delta^\Op_k n^{-1}, 
\end{equation*}
then it follows that on $\Omega_q$ we have
\begin{equation*}
\pen(k)\leq \hpen(k)\leq 30 \pen(k) \qquad \forall\;1\leq k\leq \Cy_n.
\end{equation*}
On $\Omega_{qp}=\Omega_q\cap\Omega_{p}$, we have  $\whk\leq \Cy_n$. Thus,
\begin{multline}
    \Bigl( \pen(k\vee \whk) + \hpen(k)-\hpen(\whk)
    \Bigr)\I{\Omega_{qp}}\\ \leq \Bigl( \pen(k)+\pen(\whk) +
    \hpen(k)-\hpen(\whk) \Bigr)\I{\Omega_{qp}}\\\label{D:eq:9}
     \leq 31\pen(k)\qquad \forall 1\leq k \leq \Cy_n.
\end{multline}
Furthermore, we obviously have $\widehat\Delta_k\leq n\Delta_k^\Op$
for every $1\leq k\leq \Cy_n$, which implies $\hdelta_k\leq n\,(1+\log
n)\,\delta_k^\Op$.
Consequently, $\hpen(k)\leq 540\,
\Ex[Y^2]\,n\,(1+\log n)$, because $\delta_k^\Op/n\leq
\Opd\zeta_\Opd\delta_k/n\leq \Opd\zeta_\Opd$ for all $1\leq
k\leq\Cy_n$ by~\eqref{D:eq:7} and the definition of~$\Cy_n$. On
$\Omega_q^c\cap\Omega_{p}$, we have $\whk\leq \Cy_n$ and hence
$\pen(k\vee\hk)\leq \pen(N_n)\leq 54\,\Ex[Y^2]$, which implies
\begin{equation*}\label{D:eq:10}
 ( \pen(k\vee\hk) + \hpen(k) - \hpen(\hk) ) \I{\Omega_q^c \cap \Omega_p}
 \leq 594\,\Ex[Y^2]\, n\, (1+\log n) \I{\Omega_q^c \cap \Omega_p}.
\end{equation*}
We note further that for all $\sol\in\cF_\solw^\solr$ with $
\sum_{j\in\N}\solw_j^{-1}=\Gamma <\infty$ and for all $z\in[0,1]$ we
have $|\sol(z)|^2\leq \solr \sum_{j\in\N}\solw_j^{-1}\basTr_j^2(z)\leq
2\solr\Gamma$ using the Cauchy-Schwarz inequality. 
Thereby, given $m\geq 1$ such that $\Ex[U^{2m}|W]\leq \sigma^{2m}$, it
follows  that 
\begin{equation}
  \label{D:eq:23}
\Ex [Y^{2m}|W]\leq 2^{2m}(2\solr\Gamma+\sigma^{2})^m   \mbox{ and, hence }  \Ex [Y^{2m}]\leq 2^{2m}(2\solr\Gamma+\sigma^{2})^m. 
\end{equation}

\bigskip
\noindent Now consider the decomposition
\begin{align*}
  \Ex\normV{\hsol_{\whk}-\sol}_\hw^2=
  \Ex\normV{\hsol_{\whk}-\sol}_\hw^2\I{\Omega_{qp}} +
  \Ex\normV{\hsol_{\whk}-\sol}_\hw^2\I{\Omega_{qp}^c}.
\end{align*}
It is now sufficient to show that for all $n\geq1$ and all $1\leq k\leq \Cy^l_n$ we have
\begin{gather}   \label{D:pr:th2:e1:1}
  \begin{split}
    \Ex\normV{\hsol_{\whk}-\sol}_\hw^2\I{\Omega_{qp}}\leq C\bigg\{
    \normV{\sol-\sol_k}_\hw^2 +\pen(k)+ \Opd \solr \max_{j\geq 1} \bigg[ \frac{\hw_j}{\solw_j}\min\big(1,\frac{1}{n\Opw_j}\big)\bigg]
\\\hfill  + \frac{(2\solr\Gamma+\sigma^2)^4}{n} +        \frac{(2\solr\Gamma+\sigma^2+1) \Opd\,\zeta_\Opd}{n}\,
\,  \Sigma\left(\frac{(2\solr\Gamma+\sigma^2)\zeta_\Opd+\vuz}{\vuz^2}\right)
 \bigg\},    
  \end{split}
  \\ \label{D:pr:th2:e1:3}
  \begin{split}
   \Ex\normV{\hsol_{\whk}-\sol}_\hw^2\I{\Omega_{qp}^c}  \leq
    \frac{C}{n}\,(2\solr\Gamma+\sigma^2),
  \end{split}
\end{gather}
because the result then follows using \eqref{D:eq:23}, that is,
$\pen(k)\leq 54\,(2\solr\Gamma+\sigma^2)\,\Opd\,\zeta_\Opd\,\delta_k
n^{-1}$, and by employing the monotonicity of $\hw/\solw$, that is
$\normV{\sol-\sol_{k}}_\hw^2\leq \solr \hw_k/\solw_k$.

\noindent Consider first  \eqref{D:pr:th2:e1:1}. Defining $\pen(k):=
54\,\Ex[Y^2]\,\delta^\Op_k n^{-1}$ and using the
estimate~\eqref{D:eq:6}, we have
\begin{multline*}
    \frac{1}{3} \normV{\hsol_{\hk} -\sol}^2_\hw \leq
    \frac{5}{3}\normV{\sol-\sol_k}^2_\hw +  9\bigg( \supt{k\vee \hk}
    |\skalarV{t,\tPhi_\nu}_\hw|^2 - 6\,\frac{\Ex[Y^2]\,\delta_{k\vee \hk}^\Op}{n} \bigg)_+ \\
\hfill + 
 \pen(k\vee \hk)
    + \hpen(k) - \hpen(\hk)\hfill\\[.75em]\hfill
 +      12\supt{k\vee\hk}|\skalarV{t,\hPhi_{\nu} -
    \tPhi_{\nu}}_\hw|^2\,\I{\Omega_q^c}
 +  12\supt{k\vee\hk}|\skalarV{t,\hPhi_{\mu}+\hPhi_{\gf} -
    \tPhi_{\gf}}_\hw|^2
\end{multline*}
and, hence using that $\hk\leq \Cy_n$ on $\Omega_p$ we obtain  for all $1\leq k\leq \Cy^l_n$
\begin{multline*}
\frac{1}{3} \normV{\hsol_{\whk}-\sol}_\hw^2\I{\Omega_{qp}}\leq \frac{5}{3} \normV{\sol-\sol_{k}}_\hw^2 + 9\sum_{k=1}^{\Cy_n}
    \bigg( \supt{k} |\skalarV{t,\tPhi_\nu}_\hw|^2
    - 6\frac{\Ex[Y^2]\delta^\Op_{k}}{n} \bigg)_+\hfill\\
\hfill +      12\supt{\Cy_n}|\skalarV{t,\hPhi_{\mu}+\hPhi_{\gf} -
    \tPhi_{\gf}}_\hw|^2+  \Bigl( \pen(k\vee \whk) + \hpen(k)-\hpen(\whk) \Bigr)\I{\Omega_{qp}}\\
\hfill\leq  \frac{5}{3} \normV{\sol-\sol_{k}}_\hw^2 + 9\sum_{k=1}^{\Cy_n}
    \bigg( \supt{k} |\skalarV{t,\tPhi_\nu}_\hw|^2
    - 6\frac{\Ex[Y^2]\delta^\Op_{k}}{n} \bigg)_+\hfill\\
\hfill  +      12\supt{\Cy_n}|\skalarV{t,\hPhi_{\mu}+\hPhi_{\gf} -
    \tPhi_{\gf}}_\hw|^2+   31 \pen(k),
\end{multline*}
where the last inequality follows from~\eqref{D:eq:9}. The second term
is bounded by employing Lemma~\ref{D:lem:talalem}. In order to control
the third term, apply Lemmata~\ref{D:lem:mu} and~\ref{D:lem:Q}. Consequently, combining these estimates proves inequality~\eqref{D:pr:th2:e1:1}.

\bigskip
{\noindent Consider now \eqref{D:pr:th2:e1:3}.} Let $\breve\sol_{k}:=\sum_{j=1}^k \fou{\sol}_j \1\{\hfou{\Op}_{jj}^2\geq 1/n\}\basTr_j$. It is easy to
see that $\normV{\hsol_{k}-\breve\sol_{k}}^2 \leq
\normV{\hsol_{k'}-\breve\sol_{k'}}^2$ for all $ k'\leq k$ and
$\normV{\breve\sol_{k}-\sol}^2\leq \normV{\sol}^2$ for all $k\geq 1$. Thus, using that $1\leq \whk\leq  {\Cy_n^u}$, we
can write
\begin{align*}
\Ex\normV{\hsol_{\whk}-\sol}_\hw^2\I{\Omega_{qp}^c} &\leq
2\{\Ex\normV{\hsol_{\whk}-\breve\sol_{\whk}}_\hw^2\I{\Omega_{qp}^c} +
\Ex\normV{\breve\sol_{\whk}-\sol}_\hw^2\I{\Omega_{qp}^c}\}\\
&\leq 2\bigg\{ \Ex \normV{\hsol_{\Cy_n^u}-\breve\sol_{\Cy_n^u}}_\hw^2\I{\Omega_{qp}^c} +
\normV{\sol}_\hw^2\, \P[\Omega_{qp}^c] \bigg\}.
\end{align*}
 Moreover, since  $\sup_{j\geq 1}\Ex[Y^4\basTr_j^4(W)]\leq 64
 (2\solr\Gamma+\sigma^2)^2 $ and $\Ex \basTr_j^4(W)\basTr_j^4(Z)\leq
 16 $ due to \eqref{D:eq:23}, it follows from
 Theorem~2.10 in \cite{Pet:95} that 
\begin{align*}
  \Ex \normV{\hsol_{\Cy_n^u}&-\breve\sol_{\Cy_n^u}}_\hw^2\I{\Omega_{qp}^c} \\
  &\leq 2n\sum_{j=1}^{\Cy_n^u}\hw_j\Bigl\{\Ex(\hfou{\gf}_j
  -\fou{\Op}_{jj} \fou{\sol}_j)^2\I{\Omega_{qp}^c} + \Ex(\fou{\Op}_{jj}
  \fou{\sol}_j-\hfou{\Op}_{jj} \fou{\sol}_j)^2\I{\Omega_{qp}^c}\Bigr\}\\
  &\leq 2n \Bigl\{ \sum_{j=1}^{\Cy_n^u} \hw_j \Bigl[\Ex
  \left(\hfou{\gf}_j - \fou{\gf}_j\right)^{4}\Bigr]^{1/2}
  \P[\Omega_{qp}^c]^{1/2}\\
  &\hspace{9em}  + \sum_{j=1}^{\Cy_n^u} \hw_j |\fou{\sol}_{j}|^2[\Ex
  (\hfou{\Op}_{jj}-\fou{\Op}_{jj})^4]^{1/2}\P[\Omega_{qp}^c]^{1/2}
  \Bigr\}  \\
  &\leq C n \Bigl\{ n\;(2\solr\Gamma+\sigma^2)
   + (n^{-1} \normV{\sol}_\hw^2) \Bigr\}
  \,\P[\Omega_{qp}^c]^{1/2},
\end{align*}
where we have used that $ \sum_{j=1}^{\Cy_n^u} \hw_j\leq n (\max_{1\leq j
  \leq {\Cy_n^u}}\hw_j)\leq n^2$ due to
Definition~\ref{D:def:unknown}~(ii). Since $(\hw/\solw)$ is
non-increasing, \eqref{D:pr:th2:e1:3} follows from
Lemmas~\ref{D:lem:POmegaqc} and~\ref{D:tech:res2}, which completes the proof.\qed


\subsection*{Illustration (continued)}
\label{sec:illustr-cont}

\proofof{Proposition~\ref{D:prop:ex-cont-ada-ukn}}
In the light of the proof of Proposition~\ref{D:coro:ex:lower} we
apply Theorem~\ref{D:thm:upper:unknown}, where in both cases the
additional conditions are easily verified
(Remark~\ref{D:sec:adapt-estim-under-1}) and the result follows by
an evaluation of the upper bound. Note further that $(\hw/\Opw)$
is in both cases non decreasing, and hence the second term in the upper bound of Theorem~\ref{D:thm:upper:unknown} is always smaller than the first one.\\
In case \textbf{[fs]} we have $\Cy_n^l\sim (n/ (\log
n))^{1/(2a+2s+1)}$ and $\kstar := n^{1/(2a+2p+1)}$. Note that
$\kstar\lesssim \Cy_n^l$. Thus, the upper bound is of order
$O((\kstar)^{-2(p-s)} + n^{-1})$, which equals $O(n^{-2(p-s)/(2a+2p+1)})$.\\
In case \textbf{[is]} we have \[\Cy_n^l\sim \{\log (n/ (\log
n)^{(2p+2a+1)/(2a)})\}^{1/(2a)}=(\log
n)^{1/(2a)}(1+o(1))\sim\kstar.\] Thereby, the upper bound is of order
$O((\kstar)^{-2(p-s)} + n^{-1})$, which equals $O((\log n)^{-(p-s)/a})$.\qed


\renewcommand{\tPhi}{\Phi}

\newpage
\section{Auxiliary results}
\label{D:sec:auxiliary-results}

\begin{lemma}\label{D:app:gen:upper:l2} Suppose that the distribution
  $P_U$ of $U$ belongs to $\cU_\sigma$, $\sigma>0$ and that the joint
  distribution of $(Z,W)$ satisfies Assumption~\ref{D:ass:A2}.  If in
  addition $\sol\in \cF_\bw^\br$ with $\Gamma= \sum_{j=1}^\infty
  \bw_j^{-1}<\infty $, then there exists a constant $C>0$ such that
  for all $k\in\N$ and for all $z\in\R^k$
  \begin{gather}\label{D:app:gen:upper:l2:e1:1}
    \Ex|z^t \,[B]_{\uk} |^{2} \leqslant (1/n)\, \normV{z}^2\,
    \sigma^2,\\\label{D:app:gen:upper:l2:e1:2} \Ex|z^t
    \,[S]_{\uk}|^{2}\leqslant (1/n)\, \normV{z}^2\,\eta^2\,\Gamma\,
    \normV{ \sol - \sol_k}_\bw^2 \\\label{D:app:gen:upper:l2:e2:1}
    \Ex\normV{[B]_{\uk}}^{4}\leqslant C\cdot
    \Bigl((k/n)\cdot\sigma^{2} \cdot
    \eta^2\Bigr)^2,\\\label{D:app:gen:upper:l2:e2:2}
    \Ex\normV{[S]_{\uk}}^{4}\leqslant C\cdot \Bigl((k/n)\cdot\eta^2
    \cdot \Gamma\cdot \normV{\sol-\sol_k}_\bw^2\Bigr)^2,
    \\\label{D:app:gen:upper:l2:e3} \Ex\normV{[\Xi]_{\uk}}^8\leqslant C
    \cdot \Bigl((k^2/n )\cdot \eta^{2}\Bigr)^4.
  \end{gather}
  Moreover, given a $(k\times k)$ matrix $M$, we have
  \begin{gather}\label{D:app:gen:upper:l2:e4}
    \Ex\normV{ M \{[B]_{\uk} + [S]_{\uk} \}}^{2} \leqslant
    (2/n)\,\tr(M^tM) \{\sigma^2+\eta^2\,\Gamma\, \normV{ \sol -
      \sol_k}_\bw^2\}.
  \end{gather}
\end{lemma}
\proof The proof of \eqref{D:app:gen:upper:l2:e1:1} --
\eqref{D:app:gen:upper:l2:e3} can be found in \cite{BJ:09}
and we omit the details.
The estimate \eqref{D:app:gen:upper:l2:e4} follows by applying
\eqref{D:app:gen:upper:l2:e1:1} and~\eqref{D:app:gen:upper:l2:e1:2} to
the identity $\normV{M\{\fou{B}_{\uk}+\fou{S}_{\uk}\}}^2=\sum_{j=1}^k
\normV{M_j^t\{\fou{B}_{\uk}+\fou{S}_{\uk}\}}^2$, where~$M_j$ denotes
the $j$-th column of $M^t$, which completes the
proof. \qed
\begin{lemma}\label{D:app:gen:upper:l3}
  Let $g=T\sol$ and for each $k\in\N$ denote
  $\sol_k:=[T]_{\uk}^{-1}[g]_{\uk}$. Given sequences $\Opw$ and
  $\solw$ satisfying Assumption~\ref{D:ass:minreg} let $T\in\cTdDw$ and
  $\sol\in \cF_\bw^\br$. For each strictly positive sequence
  $\hw:=(\hw_j)_{j\in\N}$ such that $\hw/\solw$ is non increasing we obtain
  for all $k\in\N$
  \begin{gather} \label{D:app:gen:upper:l3:e1} \normV{\sol-\sol_k}_\hw^2
    \leq 4\,\OpD\, \Opd \, \solr\, \frac{\hw_k}{\solw_k}\max\bigg( 1,
    \frac{\Opw_k}{\hw_k} \max_{1\leq j\leq k}\frac{\hw_j}{\Opw_j}\bigg)
  \end{gather}
\end{lemma}
\proof The condition $\Op\in\cTdDw$, that is,
$\sup_{k\in\N}\normV{[\diag(\Opw)]^{1/2}_{\uk}
  [\Op]_{\uk}^{-1}}^2\leq \OpD$ and $ \normV{ Tf}^2 \leqslant
\Opd \normV{ f}_\Opw^2$ for all $f\in L^2_Z$, together with the
identity \[\fou{E_k \sol-\sol_k}_{\uk} = -\fou{\Op}_{\uk}^{-1}\fou{\Op
  E_k^\perp \sol}_{\uk}\] imply \[\normV{E_k\sol-\sol_k}^2_\Opw
\leq \OpD \normV{ \Op E_k^\perp \sol}^2\leq \OpD \Opd \normV{
  E_k^\perp \sol}^2_\Opw\leq \OpD \Opd \solw_k^{-1}\Opw_k \solr \]for
all $\sol\in\cF_\solw^\solr$ because $(\Opw/\solw)$ is monotonically
non increasing. From this estimate we conclude
\begin{multline}\label{D:eq:20}
  \normV{E_k\sol-\sol_k}^2_{w}= \normV{[\diag(w)]_{\uk}^{1/2}\fou{E_k
      \sol-\sol_k}_{\uk}}^2\\\hfill\leq
  \normV{[\diag(w)]_{\uk}^{1/2}[\diag(\Opw)]_{\uk}^{-1/2}}^2
  \normV{E_k\sol-\sol_k}^2_\Opw \leq \OpD \Opd \solr
  \frac{\Opw_k}{\solw_k} \max_{1\leq j\leq k}\frac{\hw_j}{\Opw_j}.
\end{multline}
Furthermore, since $(\hw/\solw)$ is non increasing, we have
$\normV{E_k\sol-\sol}^2_{w}\leq \solr \hw_k/\solw_k$ for all
$f\in\cF_\solw^\solr$. The assertion follows now by combination of the
last estimate and~\eqref{D:eq:20} via a decomposition based on an
elementary triangular inequality.\qed

The next assertion is an immediate consequence of Lemma~A.3 from
\cite{BJ:09}. The proof, which is based on Bernstein's inequality, can
be found there.
\begin{lemma}\label{D:app:gen:upper:l4} Suppose that the joint
  distribution of $(Z,W)$ satisfies Assumption~\ref{D:ass:A2}. If in
  addition the sequence $\tw$ fulfills Assumption~\ref{D:ass:minreg},
  then for all $k\in\N$ we have
  \begin{gather*}
    P\left(\normV{[\Xi]_{\uk}}^2> \frac{\Opw_{k}}{4\OpD}\right)\leq 2
    \exp\left\{ -\frac{n\Opw_{k}}{32k^2 \OpD \eta^2}+ 2 \log k\right\}.
  \end{gather*}
\end{lemma}

\bigskip

\noindent The next Lemma corresponds to Lemma~A.4 in \cite{BJ:09}. We give the
details of the proof in the notation of the present paper for
convenience. 
\begin{lemma}\label{lem:a4fromcb}
Under the conditions of Theorem~\ref{D:res:upper:A3} we have for all
$n\leq 1$ that
\begin{align}\label{eq:a4a}
  (\kstar)^{12}\;P(\Omega_{1/2}^c) &\leq C(\solw,\Opw,\eta,\OpD)\\\label{eq:a4b}
  (R_n^\ast)^{-1}\;P(\Omega^c) &\leq C(\solw,\Opw,\eta,D)
\end{align}
\end{lemma}
\proof Proof of~\eqref{eq:a4a}. Since
$\normV{[\Op]^{-1}_{\uk}}^2\leq\OpD\Opw_k^{-1}$ due to $T\in\cTdDw$,
it follows from Lemma~\ref{D:app:gen:upper:l4} that
\[P(\Omega_{1/2}^c) \leq 2\exp\bigg\{ - \frac{n\Opw_{\kstar}}{32\OpD
  (\kstar)^2\eta^2} + 2 \log \kstar \bigg\}.\]
Due to condition~\eqref{eq:thm:upp} there is a $n_0\in\NN$ such that 
$n\Opw_{\kstar}\geq 448 \OpD\eta^2(\kstar)^2\log\kstar$ for all
$n\geq n_0$. Consequently, $(\kstar)^{12}P(\Omega_{1/2}^c)\leq 2$ for
all $n\geq n_0$, while trivially $(\kstar)^{12}P(\Omega_{1/2}^c)\leq
(k_{n_0}^\ast)^{12}$ for all $n<n_0$, which implies~\eqref{eq:a4a}
since $n_0$ and $k_{n_0}^\ast$ depend on $\solw, \Opw, \eta$, and
$\OpD$ only. 

\bigskip
\noindent Proof of~\eqref{eq:a4b}. Let $n_0\in\NN$ such that
$(\kstar)^2\max\{|\log\dstar|,(\log\kstar)\}\leq n\Opw_{\kstar}
(96\OpD\eta^2)^{-1}$ for all $n\geq n_0$. Observe that
$\Omega_{1/2}\subset\Omega$ if $n\geq 4 \OpD\Opw_{\kstar}^{-1}$. Since
$(\kstar)^{-2} n \Opw_{\kstar} \geq 96 \OpD\eta^2$ for all $n\geq n_0$
it follows $n\Opw_{\kstar} \geq 4\OpD$ for all $n\geq n_0$ and hence
$(\dstar)^{-1} P(\Omega^c)\leq(\dstar)^{-1}P(\Omega_{1/2}^c) \leq 2$
for all $n\geq n_0$ as in the proof of~\eqref{eq:a4a}. Combining the
last estimate and the elementary inequality $(\dstar)^{-1} P(\Omega^c)
\leq (\dstar)^{-1}$ for all $n\leq n_0$ shows~\eqref{eq:a4b} since
$n_0$ depends on $\solw, \Opw, \eta$, and $\OpD$ only. 
\qed


\renewcommand{\an}{{n^{1/3}}}

\noindent The key argument used in the proof of the next Lemma is the following
inequality due to \cite{Tal:96} (see also \cite{KR:05}, for example).
\begin{theorem}[\cite{Tal:96}]\label{thm:talagrand}
  Let $T_1, \ldots, T_n$ be independent random variables and $\nu^*_n(r)
  = (1/n)\sum_{i=1}^n\big[r(T_i) - \E[r(T_i)] \big]$, for $r$
  belonging to a countable class $\R$ of measurable functions. Then,
  \begin{align*}
    \E[\sup_{r\in\R} |\nu^*_n(r)|^2 - 6H_2^2]_+ \leq
    C\left(\frac{v}{n}\exp(-(nH_2^2/6v)) +
      \frac{H_1^2}{n^2}
      \exp(-K_2(nH_2/H_1)) \right)
  \end{align*}
  with numerical constants $K_2= (\sqrt{2}-1)/(21\sqrt{2})$ and $C>0$ and where
  \[ \sup_{r\in\R}\normV{r}_\infty \leq H_1,\quad
   \E\left[\sup_{r\in\R}|\nu^*_n(r)|\right]\leq H_2, \quad 
 \sup_{r\in\R}  \frac{1}{n}\sum_{i=1}^n \var(r(T_i))\leq v.\]
\end{theorem}

\begin{lemma}\label{D:lem:talalem}
There exists a numerical
constant $C>0$ such that
\begin{multline*}
   \sum_{k=1}^{\Cy_n} \E \Bigl[
    \bigg( \supt{k} |\skalarV{t,\Phi_\nu}_\hw|^2
    - \frac{6\,\Ex[Y^2]\,\delta_k^T}{n} \bigg)_+
    \Bigr] 
\\ \leq 
\frac{C}{n}
\Bigg\{ (2\solr\Gamma+\sigma^2+1) \Opd\,\zeta_\Opd\,
\,  \Sigma\left(\frac{(2\solr\Gamma+\sigma^2)\zeta_\Opd+\vuz}{\vuz^2}\right)
\Bigg\}.
\end{multline*}
where  $\Sigma(\cdot)$ is the function from  Definition~\ref{D:def:known}
\end{lemma}

\proof For $t\in\cS_k$, define $r_t(y,w) := \sum_{j=1}^k
\hw_j y\I{[|y|\leq\an]}\basTr_j(w)\fou{t}_j{\fou{\Op}}_{jj}^{-1}$. Then
it is readily seen that $\skalarV{t,\Phi_\nu}_\hw =
\frac{1}{n}\sum_{k=1}^n r_t(Y_k,W_k) - \E[r_t(Y_k,W_k)]$. 

Next, we compute constants $H_1$, $H_2$, and $v$ verifying the three
inequalities required in Talagrand's inequality
(Theorem~\ref{thm:talagrand}). Consider $H_1$ first:
\begin{align*}
    \supt{k}\normV{r_t}^2_\infty 
=  \sup_{y,w}  \sum_{j=1}^k \hw_j\left(  y
  \I{[|y|\leq\an]}\fou{\Op}_{jj}^{-1}\basTr_j(w)    \right)^2
\leq 2 n^{2/3}\delta^T_k =: H_1^2
  \end{align*}

Next, find  $H_2$. Notice that  
\begin{align*}
    \E[\supt{k} |\skalarV{t,\tPhi_\nu}_\hw|^2]
&= \frac{1}{n}\sum_{j=1}^k \hw_j|\fou{\Op}_{jj}|^{-2}\,
\var(Y\I{[|Y|\leq\an]}\basTr_j(W))\\
& \leq\frac{1}{n} \sum_{j=1}^k \hw_j|\fou{\Op}_{jj}|^{-2}\;
\Ex[\Ex[Y^2|W]\,\basTr_j(W)^2]
\leq 2 \Ex[Y^2]\frac{ \delta^T_k}{n}   =: H^2_2
  \end{align*}
As for $v$, we note that due to \eqref{D:eq:23} for all $\sol\in\cF_\solw^\solr$  the condition $P_U\in\cU_\sigma$, i.e., $\Ex[U^2|W]\leq \sigma^2$, implies $\Ex[Y^2|W]\leq 2(2\solr\Gamma +\sigma^2)$, and hence
\begin{align*}
  \supt{k} \var(&r_t(Y,W))
\leq \supt{k} \Ex\bigg[ \bigg( Y\sum_{j=1}^k
\frac{\hw_j\fou{t}_j}{\fou{\Op}_{jj}}\basTr_j(W)   \bigg)^2\bigg]\\
&= \supt{k} \Ex\bigg[ \E[Y^2|W]\; \bigg( \sum_{j=1}^k
\frac{\hw_j\fou{t}_j}{\fou{\Op}_{jj}}\basTr_j(W)   \bigg)^2\bigg]\\
&\leq 2(2\solr\Gamma +\sigma^2) \supt{k}\sum_{j,j'=1}^k
\frac{\hw_j\hw_{j'}\fou{t}_j\fou{t}_{j'}}{\fou{\Op}_{jj}\fou{\Op}_{j'j'}} \E[\basTr_j(W)\basTr_{j'}(W)]\\
&\leq 2(2\solr\Gamma +\sigma^2) \max_{1\leq j\leq k}\frac{\hw_j}{\fou{\Op}_{jj}^2}  \supt{k}\sum_{j=1}^k
\hw_j\fou{t}_j^2
\leq 2(2\solr\Gamma +\sigma^2) \Delta^T_k =: v,
\end{align*}
Employing  Theorem~\ref{thm:talagrand}, we conclude
\begin{multline*}
  \sum_{k=1}^{\Cy_n} \E \Bigl[
    \bigg( \supt{k} |\skalarV{t,\Phi_\nu}_\hw|^2
    - \frac{6\,\Ex[Y^2]\,\delta_k^T}{n} \bigg)_+
    \Bigr] \\
\leq 
C \,
\Bigg\{\frac{\Ex[Y^2]}{n}
\,  \sum_{k=1}^{\Cy_n} \frac{(2\solr\Gamma+\sigma^2)}{\Ex[Y^2]}\,\Delta_k^T\exp\left(-\frac{\Ex[Y^2]}{6(2\solr\Gamma+\sigma^2)}(\delta_k^T /
  \Delta_k^T) \right) 
\\\hfill+  n^{2/3}
\,\exp\left(
-K_2\,\sqrt{\Ex[Y^2]}n^{1/6}
\right)\sum_{k=1}^{\Cy_n} \frac{\delta_k^T}{n^2}
\Bigg\}.
\end{multline*}
The definition of $\Cy_n$ together with \eqref{D:eq:7} implies
$\sum_{k=1}^{\Cy_n}\delta_k^\Op/n^2\leq \zeta_\Opd$. Thereby,
 using \eqref{D:eq:7}, $\Delta_k^\Op\leq \Opd\tau_k$ and  the function $\Sigma$ given in Definition~\ref{D:def:known}, there exists a numerical constant $C>0$ such that 
\begin{multline*}\label{D:eq:15}
  \sum_{k=1}^{\Cy_n} \E \Bigl[
    \bigg( \supt{k} |\skalarV{t,\Phi_\nu}_\hw|^2
    - \frac{6\,\Ex[Y^2]\,\delta_k^T}{n} \bigg)_+
    \Bigr] \\
\hfill\leq 
\frac{C}{n}
\Bigg\{ \Ex[Y^2] \Opd
\,  \Sigma\Bigl(\frac{(2\solr\Gamma+\sigma^2)\zeta_\Opd}{\Ex[Y^2]}\Bigr)
+
\,\zeta_\Opd \Sigma\Bigl(\frac{1}{\sqrt{\Ex[Y^2]}}\Bigr)
\Bigg\}.
\end{multline*}
 Moreover, we have
 $\Ex[Y^2]\leq 2(2\solr\Gamma+\sigma^2)$
 and \[\inf_{\sol\in \cF_\solw^\solr} \Ex[Y^2]\geq \inf_{\sol\in L^2_Z}
 \Ex[\sol(Z)+U)^2]\geq \Ex [(U-\Ex[U|Z])^2]=\Ex[\var(U|Z)]=\vuz^2,\] which implies the result.\qed


\begin{lemma}\label{D:lem:mu}
For every $n\in\N$ we have
\[\E\bigg[\supt{\Cy_n}|\skalarV{t,\hPhi_\mu}_\hw|^2\bigg] \leq  2^9
  (2\solr\Gamma+\sigma^{2})^{4} n^{-1}.
\]
\end{lemma}
\proof Since $[\mu]_j=[\hmu]_j-\Ex [\hmu]_j$ and $\var\fou{\hmu}_j \leq n^{-1} \Ex[ Y^2\I{[|Y|>\an]}\basTr_j^2(W)]$,  it is easily seen that
\begin{align*}
  \E\bigg[\supt{\Cy_n} |\skalarV{t,\hPhi_\mu}_\hw|^2\bigg] 
\leq n\sum_{j=1}^{\Cy_n} \hw_j\var\fou{\widehat\mu}_j \leq \sum_{j=1}^{\Cy_n} \Ex\Bigg[ \bigg(\Ex[ Y^4|W]\Ex[\I{[|Y|>\an]}|W]\bigg)^{1/2}\basTr_j^2(W)\Bigg].
\end{align*}
Moreover,  we have  $\Ex [Y^{12}|W]\leq
2^{12}(2\solr\Gamma+\sigma^{2})^6$for all $\sol\in\cF_\solw^\solr$ and
$U\in\cU_\sigma$ due to \eqref{D:eq:23} with $m=6$, and hence by Markov's inequality
\[\Ex[\I{[|Y|>\an]}|W]\leq  2^{12}(2\solr\Gamma+\sigma^{2})^6 n^{-4}.\] Combining these estimates, we obtain 
 \begin{gather*}
  \E\bigg[\supt{\Cy_n} |\skalarV{t,\hPhi_\mu}_\hw|^2\bigg] \leq \sum_{j=1}^{\Cy_n} \Ex\Bigg[ 2^8(2\solr\Gamma+\sigma^{2})^{4} n^{-2}\basTr_j^2(W)\Bigg]\leq 2^9\Cy_n
  (2\solr\Gamma+\sigma^{2})^{4} n^{-2}.
\end{gather*}
The result follows now from $\Cy_n\leq n$.\qed


\begin{lemma}\label{D:lem:Q}
 There is a numerical constant $C>0$ such that for all $\sol\in\cF_\solw^\solr$ and every $k,n\in\N$
\[\E\bigg[\supt{k}|\qtII|^2\bigg] \leq 
C \Opd \solr \max_{j\geq 1} \bigg\{ \frac{\hw_j}{\solw_j}\min\big(1,\frac{1}{n\Opw_j}\big)\bigg\}.
\]
\end{lemma}
\proof Firstly, as $\sol\in\cF_\solw^\solr$, it is easily seen that
\begin{align*}
  \E\bigg[\supt{k} |\qtII|^2\bigg] 
\leq
\sum_{j=1}^k\fou{\sol}_j^2\hw_j\E[R_j^2]
 \leq
 \solr  \max_{j\geq 1} \bigg\{ \frac{\hw_j}{\solw_j}\E[R_j^2]\bigg\}
\end{align*}
where $R_j$ is defined by 
\begin{equation}
R_j
:=\left(\frac{\fou{\Op}_{jj}}{\hfou{\Op}_{jj}}\I{[\hfou{\Op}_{jj}^2\geq
    1/n]} -1\right).\label{D:eq:RRR} 
\end{equation}
The result follows from $\Ex R_j^2 \leq C\Opd \min\big(1,\frac{1}{n\Opw_j}\big)$, which can be shown as follows.  Consider the identity
\begin{align}\label{D:eq:18}
  \begin{split}
    \E|R_j|^2  
    = \E\bigg[\;\bigg|\frac{\fou{\Op}_{jj}}{\hfou{\Op}_{jj}}-1
    \bigg|^2\I{[\hfou{\Op}_{jj}^2\geq 1/n]} \;\bigg]
    +\P[\hfou{\Op}_{jj}^2<1/n] =: R^{I}_j + R^{II}_j.
  \end{split}
\end{align}
Trivially, $R_j^{II}\leq 1$. 
If $1\leq 4/(n\,\fou{\Op}_{jj}^2)$, then obviously $R^{II}_j\leq
 {4}/{(n\fou{\Op}_{jj}^2)}\leq 4\Opd/(n\Opw_j)$.
Otherwise, we have $1/n<\fou{\Op}_{jj}^2/4$ and hence, using \tschebby's
inequality, 
\begin{align*}
  R_j^{II} 
\leq \P[|\hfou{\Op}_{jj}-\fou{\Op}_{jj}| > |\fou{\Op}_{jj}|\,/2\,]
\leq \frac{4\,\var(\hfou{\Op}_{jj})}{\fou{\Op}_{jj}^2}
\leq \frac{16}{n\fou{\Op}_{jj}^2}\leq \frac{16\Opd}{n\Opw_j},
\end{align*}
where we have used that $\var(\hfou{\Op}_{jj})\leq 4/n$ for all $j$. Combining both estimates we have $R_j^I\leq 16 \Opd \min\big(1,\frac{1}{n\Opw_j}\big)$. 
Now consider $R^I_j$. We find that 
\begin{align}\label{D:eq:13}
  R^I_j = \E\bigg[ \frac{|\hfou{\Op}_{jj}-\fou{\Op}_{jj}|^2}{\hfou{\Op}_{jj}^2} \;
  \;\I{[\hfou{\Op}_{jj}^2\geq 1/n]} \bigg]
 \leq n\var( \hfou{\Op}_{jj})
\leq  4.
\end{align}
Using that $\E[|\hfou{\Op}_{jj}-\fou{\Op}_{jj}|^4]\leq c/n^2$ for some
 numerical constant $c>0$ (cf.~Theorem~2.10 in \cite{Pet:95}), there exists  a numerical constant $c>0$ such that
\begin{align*}
 R^I_j  &\leq \E\bigg[ \frac{|\hfou{\Op}_{jj}-\fou{\Op}_{jj}|^2}{\hfou{\Op}_{jj}^2} \;
  \;\I{[\hfou{\Op}_{jj}^2\geq1/n]} \; 2
  \bigg\{ \frac{|\hfou{\Op}_{jj} - \fou{\Op}_{jj}|^2}{\fou{\Op}_{jj}^2} +
  \frac{\hfou{\Op}_{jj}^2}{\fou{\Op}_{jj}^2} \bigg\} \bigg]\\&\leq
\frac{2\,n\,\E[|\hfou{\Op}_{jj}-\fou{\Op}_{jj}|^4]}{\fou{\Op}_{jj}^2}
+
\frac{2\; \var(\hfou{\Op}_{jj})}{\fou{\Op}_{jj}^2}
\leq
\frac{c}{n\,\fou{\Op}_{jj}^2}\leq \frac{c\Opd}{n\Opw_j}.
\end{align*}
 Combining with \eqref{D:eq:13} gives $R^I_j\leq C\Opd\min\Big\{1,
 \frac{1}{n\Opw_j}\Big\}$ for some numerical constant $C>0$, which completes the proof.\qed

\begin{lemma}\label{D:lem:Q2} There is a numerical constant $C>0$ such
  that
  \[\E\bigg[\supt{\Cy_n}|\skalarV{t,\hPhi_\nu - \Phi_\nu}_\hw\I{\Omega_q^c}|^2\bigg]
{\leq Cd (\P[\Omega_q^c])^{(1/2)}}
.\]
\end{lemma}
\proof  
Given $R_j$ from~\eqref{D:eq:RRR}  we begin our proof observing that
\begin{align*}
  \E\bigg[\supt{\Ce_m}|\skalarV{t,\hPhi_\nu - \Phi_\nu}_\hw\I{\Omega_q^c}|^2\bigg]
  &\leq \sum_{j=1}^{\Cy_n} \frac{\hw_j}{\fou{\Op}_{jj}^2}
  \;\E[\fou{\nu}_j^2\,R_j^2\,\I{\Omega_q^c}]\\
&\leq \sum_{j=1}^{\Cy_n} \frac{\hw_j}{\fou{\Op}_{jj}^2}
  \;\big(\E[\fou{\nu}_j^8]\E[R_j^8]\big)^{1/4}\; \P[{\Omega_q^c}]^{1/2},
\end{align*}
where we have applied Cauchy-Schwarz twice. By Petrov's
inequality, there exists a numerical constant $c>0$ such that $E[\fou{\nu}_j^8] \leq cn^{-4/3}$
and hence, because $d\delta_k \geq \sum_{j=1}^k \frac{\hw_j}{\fou{\Op}_{jj}^2}$,
\[  \E\bigg[\supt{\Ce_m}|\skalarV{t,\hPhi_\nu -
  \Phi_\nu}_\hw\I{\Omega_q^c}|^2\bigg]
 \leq  \P[\Omega_q^c]^{1/2}d\delta_k \max_{1\leq j\leq
   \Cy_n}(\E[R_j^8])^{1/4} \]
In analogy to~\eqref{D:eq:18}, we decompose the moment of $R_j$ into two terms 
\[\E[R_j^8] = \E\bigg[\;\bigg|\frac{\fou{\Op}_{jj}-\hfou{\Op}_{jj}}{\hfou{\Op}_{jj}}
    \bigg|^8\I{[\hfou{\Op}_{jj}^2\geq 1/n]} \;\bigg]
    +\P[\hfou{\Op}_{jj}^2<1/n],\]
which we bound by a constant using Petrov's inequality.  \qed

\begin{lemma}\label{D:lem:POmegaqc}  We have $\P[\Omega_q^c]\leq
  2(2016\Opd/\Opw_1)^7\, n^{- 6}$, where $\Omega_q$ is the event defined
  in~\eqref{eq:OmOmpq}.
\end{lemma}
\proof
Consider the complement of $\Omega_q$ given by 
\[\Omega_q^c = \bigg\{\exists\; 1\leq j \leq \Cy_n  \; \bigg|\;
\Big|\frac{\fou{\Op}_{jj}}{\hfou{\Op}_{jj}}-1\Big| > \frac{1}{2} \;\vee\;
 \hfou{\Op}_{jj}^2 <1/n    \bigg\}.\]
It follows from Assumption~\ref{D:ass:known} (i) that  $\fou{\Op}_{jj}^2\geq2/n$ for all
$ 1\leq j\leq\Cy_n$. This yields
\[\P(\Omega_q^c)\leq \sum_{j=1}^{\Cy_n} \P\bigg[\bigg|\frac{\hfou{\Op}_{jj}}{\fou{\Op}_{jj}} -1  \bigg|>\frac{1}{3}\bigg].\]
From Hoeffding's inequality follows
\begin{equation*}
\P[|\hfou{\Op}_{jj}/\fou{\Op}_{jj} -1
|>1/3]\leq2\,\exp\bigg(-\frac{n\fou{\Op}_{jj}^2}{288}\bigg),
 \end{equation*}
which implies the result by definition of $\Cy_n$.\qed

\begin{lemma}\label{D:tech:res2}
Consider the event $\Omega_p$ defined in~\eqref{eq:OmOmpq}. Then  we
have
 \[\P(\Omega_{p}^c) \leq 4\,\left(\frac{2016\,\Opd}{\Opw_1}\right)^7n^{-6}, \qquad\forall\; n\geq 1.\]
\end{lemma}
\proof  Let $\Omega_I:=\{ \Cy^l_n>\hCy_n\}$ and $\Omega_{II}:=\{ \hCy_n > \Cy_n\}$. Then we have $\Omega_{p}^c=\Omega_I\cup \Omega_{II}$.  
Consider $\Omega_I$ first.  By definition of $\Cy^l_n$,
we have that $\min_{1\leq j\leq\Cy^l_n}
\frac{|[\Op]_j|^{2}}{|j|(\hw_j)_{\vee1}}\geq \frac{4 (\log n)}{n}$, which implies
   \begin{multline*}
     \{ \hCy_n< \Cy^l_n\}\subset \bigg\{ \exists 1\leq
     j\leq \Cy^l_n \;\bigg|\;
     \frac{\hfou{\Op}_{jj}^2}{|j|(\hw_j)_{\vee1}}<\frac{\log n}{n}\bigg\}\\\subset
     \bigcup_{1\leq j\leq \Cy^l_n}\bigg\{
     \frac{|\hfou{\Op}_{jj}|}{|\fou{\Op}_{jj}|}\leq 1/2\bigg\}\subset
     \bigcup_{1\leq j\leq \Cy^l_n}\bigg\{
     |\hfou{\Op}_{jj}/\fou{\Op}_{jj} -1   |\geq 1/2\bigg\}.
  \end{multline*}
Therefore, $\Omega_I \subset \bigcup_{1\leq|j|\leq 
  \Cy_n}\Bigl\{ |\hfedf_j/\fedf_j-1|\geq 1/2\Bigr\}$, since $\Cy_n^l\leq \Cy_n$. Hence, as in~\eqref{D:eq:7}
applying Hoeffding's inequality together with the definition of
$\Cy_n$ gives
\begin{equation}\label{D:tech:e1}
  \P[\Omega_I] \leq 
 \sum_{j=1}^{\Cy_n} 2\,\exp\bigg(-\frac{n\,\fou{\Op}_{jj}^2}{288}\bigg)
\leq 2 \left(\frac{2016\,\Opd}{\Opw_1} \right)^7n^{-6}.\end{equation}
Consider $\Omega_{II}$.
 Recall that $\frac{\log n}{4n} \geq
 \max_{|j|\geq \Cy_n} \frac{\fou{\Op}_{jj}^2}{|j|(\hw_j)_{\vee1}}$ due to Assumption~\ref{D:ass:unknown}, and hence
 \begin{align*}
   \{\hCy_n> \Cy_n\}&\subset \Bigl\{ \forall 1\leq j\leq
   \Cy_n\;\Big|\; \frac{\hfou{\Op}_{jj}^2}{|j|(\hw_j)_{\vee1}}\geq \frac{\log n}{n}\Bigr\}\\ &\subset
   \biggl\{ \frac{|\hfou{\Op}_{\Cy_n}|}{|\fou{\Op}_{\Cy_n}|}\geq
   2\biggr\}\subset \Bigl\{ |\hfou{\Op}_{\Cy_n}/\fou{\Op}_{\Cy_n}-1|\geq
   1\Bigr\}.
  \end{align*}
  Hoeffding's inequality and the definition of $\Cy$ yield
  $ \P[\Omega_{II}] \leq 2(2016\Opd/\Opw_1)^7n^{-6}$, which by
  combining with \eqref{D:tech:e1} implies the result.  \qed


\bibliographystyle{apalike}
\bibliography{d.bib}

\end{document}